\documentclass{amsart}          

\usepackage{amsmath,amsthm}     
\usepackage{amssymb}            
\usepackage{euscript}           
\usepackage{enumerate,calc}
\usepackage[matrix,arrow,curve,frame]{xy}    

\xymatrixcolsep{1.9pc}                          
\xymatrixrowsep{1.9pc}
\newdir{ >}{{}*!/-5pt/\dir{>}}                  

\def\longfib{\DOTSB\relbar\joinrel\twoheadrightarrow}

\newtheorem{thm}[subsection]{Theorem}
\newtheorem{defn}[subsection]{Definition}
\newtheorem{prop}[subsection]{Proposition}

\newtheorem{cor}[subsection]{Corollary}
\newtheorem{lemma}[subsection]{Lemma}

\theoremstyle{definition}  
\newtheorem{example}[subsection]{Example}

\newtheorem{remark}[subsection]{Remark}

\newcommand{\dfn}{\textbf} 
\newcommand{\mdfn}[1]{\dfn{\mathversion{bold}#1}} 

\newcommand{\tens}              {\otimes}               
\newcommand{\extens}            {\boxtimes}
\newcommand{\iso}               {\cong}

\newcommand{\cat}{\EuScript}    
\newcommand{\cA}{{\cat A}}      
\newcommand{\cB}{{\cat B}}      
\newcommand{\cC}{{\cat C}}
\newcommand{\cD}{{\cat D}}
\newcommand{\cE}{{\cat E}}
\newcommand{\cF}{{\cat F}}

\newcommand{\cI}{{\cat I}}

\newcommand{\cM}{{\cat M}}
\newcommand{\cN}{{\cat N}}
\newcommand{\cO}{{\cat O}}
\newcommand{\cP}{{\cat P}}
\newcommand{\cQ}{{\cat Q}}
\newcommand{\cR}{{\cat R}}

\newcommand{\cT}{{\cat T}}
\newcommand{\cV}{{\cat V}}

\newcommand{\Set}{{\cat Set}}

\newcommand{\sSet}{s{\cat Set}}
\newcommand{\dMod}{{\cat Mod}}

\newcommand{\Ho}{{\cat Ho}}
\newcommand{\ho}{\text{Ho}\,}

\newcommand{\field}[1]  {\mathbb #1} 

\newcommand{\F}         {\field F}

\newcommand{\Z}         {\field Z}

\newcommand{\Q}         {\field Q}

\DeclareMathOperator*{\colim}{colim}

\DeclareMathOperator{\spec}{Sp^{\Sigma}}
\DeclareMathOperator{\Hom}{Hom}

\DeclareMathOperator{\vMap}{vMap}

\DeclareMathOperator{\diag}{diag}
\DeclareMathOperator{\coeq}{coeq}

\DeclareMathOperator{\Mod}{Mod-}

\newcommand{\ra}{\rightarrow}                   
\newcommand{\lra}{\longrightarrow}              
\newcommand{\lla}{\longleftarrow}               
\newcommand{\llra}[1]{\stackrel{#1}{\lra}}      
\newcommand{\llla}[1]{\stackrel{#1}{\lla}}      

\newcommand{\we}{\llra{\sim}}                   
\newcommand{\bwe}{\llla{\sim}}
\newcommand{\cof}{\rightarrowtail}              
\newcommand{\trfib}{\stackrel{\sim}{\longfib}}
\newcommand{\trcof}{\stackrel{\sim}{\cof}}

\newcommand{\inc}{\hookrightarrow}              
\newcommand{\dbra}{\rightrightarrows}           

\newcommand{\blank}{-}                          
\newcommand{\Id}{Id}                            
\newcommand{\id}{id}                            
\newcommand{\und}{\underline}
\newcommand{\mM}{\underline{\cM}}

\newcommand{\bd}{\partial}

\newcommand{\adjoint}{\rightleftarrows}
\newcommand{\ladjoint}[1]{\stackrel{#1}{\adjoint}}      

\newcommand{\bdd}[1]{\partial\Delta^{#1}}
\newcommand{\del}[1]{\Delta^{#1}}

\newcommand{\he}{\simeq}

\numberwithin{equation}{section}

\newenvironment{myequation}
  {\addtocounter{subsection}{1}\begin{equation}}
  {\end{equation}$\!\!$}

\newcommand{\DGA}{\cD G\cA}

\newcommand{\Ch}{Ch}
\newcommand{\ch}{ch}

\newcommand{\U}{U}
\newcommand{\Ua}{U_{ad}}
\newcommand{\Spe}{Sp}
\newcommand{\Spes}{Sp^\Sigma}
\DeclareMathOperator{\hEnd}{hEnd}

\newcommand{\mC}{\underline{\cC}}
\newcommand{\mN}{\underline{\cN}}
\newcommand{\unit}{\mathbf{1}}

\DeclareMathOperator{\hdga}{\hEnd_{ad}}
\DeclareMathOperator{\hdgadga}{\hEnd_{dga}}
\newcommand{\Ab}{\cat Ab}
\DeclareMathOperator{\Ad}{Ad}

\newcommand{\ME}{ME}

\newcommand{\cofib}{\cof}
\renewcommand{\F}{F}
\newcommand{\Ring}{\cat{R}ing}
\newcommand{\Fr}{Fr}
\newcommand{\ccO}{\overline{\cO}}
\newcommand{\Cat}{{\cat Cat}}
\newcommand{\Graph}{{\cat Graph}}
\newcommand{\Monoid}{{\cat Monoid}}

\newcommand{\namma}{\nut}
\newcommand{\nut}{\nu}

\newcommand{\ff}{\cF}

\DeclareMathOperator{\Fhomu}{\mM_{\Ch}}

\DeclareMathOperator{\ND}{\mN_{\cD}}
\DeclareMathOperator{\Nsab}{\mN_{\spec(\sab)}}
\newcommand{\ehom}{{\Hom}}

\newcommand{\ChZ}{\Ch}

\newcommand{\Sp}{Sp}

\newcommand{\sab}{s{\cat Ab}}
\DeclareMathOperator{\Func}{Func}
\DeclareMathOperator{\Sing}{Sing}
\DeclareMathOperator{\Real}{Re}

\begin{document}

\title[Enriched model categories]{Enriched model categories and an
application to additive endomorphism spectra}

\author{Daniel Dugger}
\address{Department of Mathematics\\ University of Oregon\\ Eugene, OR
97403} 

\email{ddugger@math.uoregon.edu}

\author{Brooke Shipley}
\address{Department of Mathematics\\ University of Illinois at Chicago
\\ Chicago, IL 60607}

\email{bshipley@math.uic.edu}

\begin{abstract}
We define the notion of an additive model category and prove that any
stable, additive, combinatorial model category $\cM$ has a model
enrichment over $\Sp^\Sigma(s\Ab)$ (symmetric spectra based on
simplicial abelian groups).  So to any object $X\in \cM$ one can
attach an endomorphism ring object, denoted $\hdga(X)$, in this
category of spectra.  One can also obtain an associated differential
graded algebra carrying the same information.  We prove that the
homotopy type of $\hdga(X)$ is an invariant of Quillen equivalences
between additive model categories.

We also develop a general notion of an adjoint pair of functors being
a `module' over another such pair; we call such things adjoint
modules.  This is used to show that one can transport enrichments over
one symmetric monoidal model category to a Quillen equivalent one, and
in particular it is used to compare enrichments over
$\Sp^\Sigma(s\Ab)$ and chain complexes.
\end{abstract}

\maketitle

\tableofcontents


\section{Introduction}
\label{se:intro}

A model category is called \dfn{additive} if two conditions are satisfied.
First, its hom-sets must have natural structures of abelian groups
with respect to which composition is biadditive.  Secondly, the
abelian group structures on these hom-sets must interact well with the
notion of `higher homotopies'.  We give a precise definition in
Section~\ref{se:additivedef}.  Examples of additive model categories
include chain complexes over a ring and differential graded modules over a 
differential graded algebra, as one should expect.

Recall that a category is \dfn{locally presentable} if it is
cocomplete and all objects are small in a certain sense; see~\cite{AR}.
A model category is called \dfn{combinatorial} if it is
cofibrantly-generated and its underlying category is locally
presentable.  
A model category is \dfn{stable} if it is pointed and
the suspension functor is an auto-equivalence of the homotopy
category.  In \cite{hend} it was shown that any stable, combinatorial
model category could be naturally enriched over the category
$\Sp^\Sigma$ of symmetric spectra.  This enrichment is invariant under
Quillen equivalences in a certain sense.

In the present paper we extend the results of \cite{hend} to show that
any stable, combinatorial, additive model category has a natural
enrichment over $\Sp^\Sigma(s\Ab)$, the category of symmetric spectra based
on simplicial abelian groups.  This enrichment is not an
invariant of Quillen equivalence, but it {\it is\/} preserved by
Quillen equivalences which only involve {\it additive\/} model
categories.

\begin{remark}
The tools developed in this paper are applied in \cite{tpwe}.  Two
additive model categories $\cM$ and $\cN$ are called \dfn{additively
Quillen equivalent} if there is a zig-zag of Quillen equivalences
between $\cM$ and $\cN$ in which every intermediate step is additive.
It is a strange fact, established in \cite{tpwe}, that additive model
categories can be Quillen equivalent but not additively Quillen
equivalent.  The demonstration of this fact uses
the model enrichments developed in the present paper.
\end{remark}

We should explain up front that there are really three
separate things going on in this paper.  One is the development of the
theory of additive model categories, taken up in
Sections~\ref{se:additivedef} and \ref{se:univ}.  The second is the
construction of the model enrichment by $\Sp^\Sigma(s\Ab)$, which is
begun in Section~\ref{se:additive}.  Most of the details of the model
enrichment exactly follow the pattern in \cite{hend}.  There is one
extra result we wish to consider, though, which involves comparing
model enrichments over $\Sp^\Sigma(s\Ab)$ to model enrichments over
the Quillen equivalent category $\Ch$ of chain complexes of abelian
groups.  For this last issue we need to develop quite a bit more about
enriched model categories than is available in the literature.  Since
this foundational material is important in its own right, we include
it at the very beginning as Sections~\ref{se:enrich} through
\ref{se:transfer}.

\subsection{A closer look at the results}  
To describe the results in more detail we need to recall some enriched
model category theory; specifically, we need the notions of \dfn{model
enrichment} and \dfn{quasi-equivalence} from \cite{hend}.  Let $\cM$
be a model category and $\cV$ be a symmetric monoidal model category.
Briefly, a model enrichment is a bifunctor $\tau\colon \cM^{op}\times
\cM \ra \cV$ together with composition maps $\tau(Y,Z)\tens \tau(X,Y)
\ra \tau(X,Z)$ which are associative and unital.  The bifunctor must
interact well with the model category structure---see \cite{hend} for
an explicit list of the necessary axioms, or Section
\ref{se:modelenrich} for a summary.

There is a notion of when two model enrichments of $\cM$ by $\cV$ are
`quasi-equivalent', which implies that they carry the same homotopical
information.  This takes longer to describe, but the reader can again
find it in Section~\ref{se:modelenrich}.  We let $\ME_0(\cM,\cV)$
denote the quasi-equivalence classes of model enrichments.

If $L\colon \cM \adjoint \cN\colon R$ is a Quillen pair, there are
induced functors denoted $L_*\colon \ME_0(\cM,\cV) \ra \ME_0(\cN,\cV)$
and $L^*\colon \ME_0(\cN,\cV) \ra \ME_0(\cM,\cV)$.  
When $(L,R)$ is a
Quillen equivalence these are inverse bijections.

Using the above language, we can state the basic results.  These are
proved in Sections~\ref{se:additive} and ~\ref{se:chain}.

\begin{thm}
\label{th:additive-enrich}
If $\cM$ is a stable, additive, combinatorial model category, then
there is a canonical element $\sigma_{\cM}\in
\ME_0(\cM,\Sp^\Sigma(s\Ab))$.  If $L\colon \cM \ra \cN$ is a Quillen
equivalence then $L_*(\sigma_\cM)=\sigma_\cN$ and $L^*(\sigma_\cN)=\sigma_\cM$.
\end{thm}

If $X\in \cM$, choose a cofibrant-fibrant object $\hat{X}$ which is
weakly equivalent to $X$.  Then $\sigma_\cM(\hat{X},\hat{X})$ gives a
ring object in $\Sp^\Sigma(s\Ab)$.  The resulting isomorphism class in
the homotopy category $\ho(\Ring[\Sp^\Sigma(s\Ab)])$ only depends on
the homotopy type of $X$.  We write $\hdga(X)$ for any ring object in
this isomorphism class, and call it the \dfn{additive homotopy
endomorphism spectrum of $X$}.

\begin{prop}
\label{pr:hdga-invariance}
Let $\cM$ and $\cN$ be additive, stable, combinatorial model
categories.  Suppose $\cM$ and $\cN$ are Quillen equivalent through a
zig-zag of additive (but not necessarily combinatorial) model
categories.  Let $X\in \cM$, and let $Y\in \ho(\cN)$ correspond to $X$
under the derived equivalence of homotopy categories.  Then $\hdga(X)$
and $\hdga(Y)$ are weakly equivalent in $\Ring[\Spe^\Sigma(s\Ab)]$.
\end{prop}

Any ring object $R$ in $\Sp^\Sigma(s\Ab)$ gives rise to a ring object
in $\Sp^\Sigma$ by forgetting the abelian group structure---this is
called the Eilenberg-Mac\,Lane spectrum associated to $R$.  Recall that
in \cite{hend} it was shown how to attach to any $X$ in a stable,
combinatorial model category an isomorphism class in
$\ho(\Ring[\Sp^\Sigma])$.  This was called the \dfn{homotopy
endomorphism spectrum of $X$}, and denoted $\hEnd(X)$.  We have the
following:

\begin{prop}
\label{pr:hdga=EM}
Given $X\in \cM$ as above, the homotopy endomorphism spectrum
$\hEnd(X)$ is the Eilenberg-Mac\,Lane spectrum associated to 
$\hdga(X)$.
\end{prop}

Finally, we have two results explaining how to compute $\hdga(X)$ when
the model category $\cM$ has some extra structure.  Recall that if
$\cC$ is a symmetric monoidal model category then a \mdfn{$\cC$-model
category} is a model category equipped with compatible tensors,
cotensors, and enriched hom-objects over $\cC$ satisfying the analogue
of $SM7$. See Section~\ref{se:enrich} for more detailed information.
For $X,Y\in \cM$ we denote the enriched hom-object by
$\und{\cM}_{\cC}(X,Y)$.

Note that a $\Spe^\Sigma(s\Ab)$-model category is automatically
additive and stable.  This follows from
Corollary~\ref{co:Cmodel=additive} below, and the appropriate
analogue of~\cite[3.5.2]{SS2} 
or~\cite[3.2]{GS}. 

\begin{prop}
\label{pr:hdga-spectral}
Let $\cM$ be a combinatorial 
$\Spe^\Sigma(s\Ab)$-model category.  Let $X\in \cM$ be
cofibrant-fibrant.  Then $\hdga(X)$ is weakly equivalent
to the enriched hom-object $\und{\cM}_{\Spe^\Sigma(s\Ab)}(X,X)$.
\end{prop}

In ~\cite{S} it is shown that the model categories of rings in
$\Spe^\Sigma(s\Ab)$ and in $\ChZ$ are Quillen equivalent.  This is
recalled in Section~\ref{se:chain}.  Note that the rings in $\ChZ$ are
just differential graded algebras (dgas).  The associated derived functors will be denoted $H'\colon
\DGA \adjoint \Ring\Spe^\Sigma(s\Ab)\colon \Theta'.$ (The reason for
the `primes' is that in \cite{S} the functors $H$ and $\Theta$ are
functors between $\DGA$ and $H\Z$-algebras with
$\Ring\Spe^\Sigma(s\Ab)$ an intermediate category.)  We then define
the \dfn{homotopy endomorphism dga of X} to be $\Theta'[\hdga(X)]$ and
write \mdfn{$\hdgadga(X)$}.  Obviously, this carries exactly the same
information as $\hdga(X)$.  In fact, $H'[\hdgadga(X)]$ is weakly
equivalent to $\hdga(X)$ since $H'$ and $\Theta'$ are inverse
equivalences on the homotopy category level.

As above, we remark that a $\Ch$-model category is automatically
additive and stable, by Corollary~\ref{co:Cmodel=additive} and
the appropriate analogue of~\cite[3.5.2]{SS2}. 

\begin{prop}
\label{pr:hdga-chain}
Let $\cM$ be a combinatorial 
$\ChZ$-model category.  Assume $\cM$ has a generating set of compact
objects, as defined in (\ref{def:CompGen}) below.  Let $X\in \cM$ be
cofibrant-fibrant.  Then $\und{\cM}_{\Ch}(X,X)$ is weakly equivalent
to $\hdgadga(X)$.
\end{prop}

The assumption about the generating set in the above proposition is
probably unnecessary, but we don't know how to remove it.  It is
satisfied in most cases of interest.

The proof of Proposition~\ref{pr:hdga-chain} is not hard, but it
requires a careful comparison of enrichments over $\Ch$ and
$\Sp^\Sigma(s\Ab)$.  This reduces to an abstract problem in enriched
model category theory, but the necessary tools do not seem to be
available in the literature.  The first part of the paper is spent
developing them.  Among other things, one needs a notion of an adjoint
pair of functors being a `module' over another such pair; we call such
things {\it adjoint modules\/}, and develop their basic theory in
Sections~\ref{se:adjmod}--\ref{se:adjmodapp}.  This notion has other
applications, most notably in \cite{GS}.

\begin{remark}
The study of \mdfn{$dg$-categories} seems to be of current
interest---see, for example, \cite{Dr, T}.  A $dg$-category is simply
a category enriched over unbounded chain complexes $\Ch_k$, where
$k$ is some commutative ground ring.  We remark that the homotopy
theory of $dg$-categories over $\Z$ is essentially the same as that of
$\Sp^\Sigma(s\Ab)$-categories
(this follows from results of \cite{S} and
\cite{SS3}).  So the present paper may be regarded as associating to
any stable, additive model category an underlying dg-category.
\end{remark}

\subsection{Organization of the paper}
Section~\ref{se:enrich} recalls the basics of enriched model category
theory as used in \cite{hend}.  
The new work begins in
Sections~\ref{se:adjmod} and \ref{se:adjmodapp} where we develop the
notion of adjoint modules.  This is used in Section~\ref{se:transfer} to
prove a technical theorem about transporting enrichments over one
symmetric monoidal model category to a Quillen equivalent one.
Sections~\ref{se:additivedef} through ~\ref{se:chain} contain the main
results on additive model categories and
$\Sp^\Sigma(s\Ab)$-enrichments. 
Appendix~\ref{se:CI-cat} reviews and expands
on material from \cite{SS3}, which is needed in 
Section~\ref{se:app-module-cat}.

\subsection{Notation and terminology}
This paper is a sequel to \cite{hend}, and we will assume the reader
is familiar with the machinery developed therein.  In particular, we
assume a familiarity with model enrichments and quasi-equivalences;
see Section~\ref{se:modelenrich} for quick summaries, though.
We use one piece of terminology which is not quite standard.  Namely,
if $\cM$ and $\cN$ are model categories then by a \dfn{Quillen map}
$L\colon \cM \ra \cN$ we mean an adjoint pair of Quillen functors
$L\colon \cM \adjoint \cN\colon R$, where $L$ is the left adjoint.


\section{Enriched model categories}
\label{se:enrich}

In this section we review the notion of a model category $\cM$ being
enriched over a second model category $\cC$.  This situation comes in
two varieties.  If for every two objects $X,Y\in \cM$ one has a
`mapping object' $\mM_{\cC}(X,Y)$ in $\cC$ together with composition
maps (subject to certain axioms), then this is called a model
enrichment.  If for every $X\in \cM$ and $c\in \cC$ one also has
objects $X\tens c$ and $\F(c,X)$ in $\cM$, related by adjunctions to
the mapping objects and also subject to certain axioms, then we say
that $\cM$ is a $\cC$-model category.  Thus, a $\cC$-model category
involves a model enrichment plus extra data.

There are two main examples to keep in mind.  A {\it simplicial\/}
model category is just another name for an $\sSet$-model category.
And if $\cM$ is any model category, then the hammock localization of
Dwyer-Kan \cite{DK} is an example of a model enrichment of $\cM$ over
$\sSet$.

\subsection{Symmetric monoidal model categories}

Let $\cC$ be a closed symmetric monoidal category.  This says that we
are given a bifunctor $\tens$, a unit object $\unit_{\cC}$, together
with associativity, commutativity, and unital isomorphisms making
certain diagrams commute (see \cite[Defs. 4.1.1, 4.1.4]{hovey} for a
nice summary).  The `closed' condition says that there is also a
bifunctor $(a,b) \mapsto \mC(a,b)\in \cC$ together with a natural
isomorphism
\[ \cC(a, \mC(b,c)) \iso \cC(a\tens b, c).
\]
Note that this gives  isomorphisms
$\cC(\unit_{\cC},\mC(a,b)) \iso \cC(\unit_{\cC}\tens a,b) \iso \cC(a,b)$.

A \dfn{symmetric monoidal model category} consists of a closed
symmetric monoidal category $\cC$, together with a model structure on
$\cC$, satisfying two conditions:
\begin{enumerate}[(1)]
\item The analogue of SM7, as given in either \cite[4.2.1]{hovey} or
\cite[4.2.2(2)]{hovey}. 
\item A unit condition given in \cite[4.2.6(2)]{hovey}.
\end{enumerate}

\subsection{$\cC$-model categories}

Let $\cC$ be a symmetric monoidal category.
One defines a \mdfn{closed $\cC$-module category} to be a category $\cM$
equipped with natural constructions which assign
to every $X,Z\in \cM$ and $c\in \cC$ objects
\[ X\tens c \in \cM, \qquad \F(c,Z) \in \cM, \qquad \text{and} \qquad 
\mM_{\cC}(X,Z) \in \cC.
\]
One requires, first, that there are natural isomorphisms $(X\tens
a)\tens b\iso X\tens (a\tens b)$ and $X\tens \unit_{\cC} \iso X$
making certain diagrams commute (see \cite[Def. 4.1.6]{hovey}).  One
of these diagrams is a pentagon for four-fold associativity.  We also
require natural isomorphisms
\begin{equation}
\label{eq:adjoint}
 \cM(X\tens a,Z) \iso \cM(X, \F(a,Z)) \iso \cC(a, \mM_{\cC}(X,Z))
\end{equation}
(see \cite[4.1.12]{hovey}).

Finally, suppose $\cC$ is a symmetric monoidal model category.  A
\mdfn{$\cC$-model category} is a model category $\cM$ which is also a
closed $\cC$-module category and where the two conditions from
\cite[4.2.18]{hovey} hold: these are again the analogue of SM7 and a unit
condition.

\subsection{Model enrichments}
\label{se:modelenrich}
Let $\cM$ be a model category and let $\cC$ be a symmetric monoidal
model category.  Recall from \cite[3.1]{hend} that a \dfn{model enrichment}
of $\cM$ by $\cC$ is a bifunctor $\sigma\colon \cM^{op} \times \cM \ra
\cC$ which is equipped with composition pairings $\sigma(Y,Z)\tens
\sigma(X,Y) \ra\sigma(X,Z)$ and unit maps $\unit_{\cC}\ra \sigma(X,X)$
satisfying associativity and unital conditions.  There is also a
compatibility condition between the functor structure and the unit
maps.  Finally, one assumes that if $X\ra X'$ is a weak equivalence
between cofibrant objects and $Y\ra Y'$ is a weak equivalence between
fibrant objects then the maps $\sigma(X,Y) \ra \sigma(X,Y')$ and
$\sigma(X',Y)\ra \sigma(X,Y)$ are weak equivalences.  See
\cite[Section 3.1]{hend}.

There is a notion of quasi-equivalence encoding when two model
enrichments are `the same'.  This is also given in \cite[Section 3.1]{hend}.
To define this we need two preliminary notions.

Let $\sigma$ and $\tau$ be two model enrichments of $\cM$ by $\cC$.
By a $\sigma-\tau$ \dfn{bimodule} we mean a collection of objects
$M(a,b)\in \cC$ for every $a,b\in \cC$, together with multiplication
maps
\[ \sigma(b,c) \tens M(a,b) \ra M(a,c) \qquad\text{and}
\qquad M(b,c)\tens \tau(a,b)
\ra M(a,c) 
\]
which are natural in $a$ and $c$.  Associativity and unital conditions
are again assumed, although we will not write these down.  One also
requires that for any $a,b,c,d\in \cC$ the two obvious
maps
\[ \sigma(c,d) \tens M(b,c) \tens \tau(a,b) \dbra M(a,d)
\]
are equal.

It is perhaps not quite obvious, but $M$ becomes a bifunctor via the
multiplication maps from $\sigma$ and $\tau$ and the fact that 
$\sigma$ and $\tau$ are bifunctors.  See \cite[Section 2.2]{hend}.

A \mdfn{pointed $\sigma-\tau$ bimodule} is a bimodule $M$
together with a collection of maps $\unit_\cC \ra M(c,c)$ for every $c\in \cC$,
such that for any map $a\ra b$ the square
\[ \xymatrix{ \unit_\cC \ar[r] \ar[d] & M(a,a) \ar[d] \\
              M(b,b) \ar[r] & M(a,b)
}
\]
commutes.

A \dfn{quasi-equivalence} between two model enrichments $\sigma$ and
$\tau$ consists of a pointed $\sigma-\tau$ bimodule $M$ such that the
compositions
\begin{align*}
\sigma(a,b)\tens \unit_\cC \ra \sigma(a,b) &\tens M(a,a) \ra M(a,b)
\qquad\text{and} \\
 \unit_\cC \tens \tau(a,b) \ra M(b,b) &\tens \tau(a,b) \ra M(a,b)
\end{align*}
are weak equivalences whenever $a$ is cofibrant and $b$ is fibrant.

The notion of quasi-equivalence generates an equivalence relation on
the class of model enrichments of $\cM$ by $\cC$.  We write
$\ME_0(\cM,\cC)$ for the collection of equivalence classes of model
enrichments.  When we say that two enrichments $\sigma$ and $\tau$ are
`quasi-equivalent' we mean that they are in the same equivalence
class; note that this means there is a chain of model enrichments
$\sigma=\sigma_1,\sigma_2,\ldots,\sigma_n=\tau$ and pointed
$\sigma_i-\sigma_{i+1}$ bimodules $M_i$ giving quasi-equivalences
between each step in the chain.

If $L\colon \cM \ra \cN$ is a Quillen map then by \cite[Prop. 3.14]{hend}
there are induced maps $L_*\colon \ME_0(\cM,\cC) \ra \ME_0(\cN,\cC)$
and $L^*\colon \ME_0(\cN,\cC) \ra \ME_0(\cM,\cC)$.
When $L$ is a Quillen equivalence these are inverse bijections.

\subsection{Monoidal functors}\label{se:monoidalfunc}

Suppose that $\cC$ and $\cD$ are symmetric monoidal model categories,
and that $F\colon \cC \adjoint \cD\colon G$ is a Quillen pair.  

First of all, recall that $G$ is called {\bf lax monoidal} if there is
a natural transformation
\[ G(X)\tens G(Y) \ra G(X\tens Y) \]
and a map $\unit_\cC \ra G(\unit_\cD)$ which are compatible with the
associativity and unital isomorphisms in $\cC$ and $\cD$.
A lax monoidal functor takes monoids in $\cD$ to monoids in $\cC$.

A lax monoidal functor is called {\bf strong monoidal} if the above
maps are actually isomorphisms.

If $G$ is lax monoidal then the adjunction gives rise to induced maps
$F(\unit_\cC) \ra \unit_\cD$ and $F(A\tens B) \ra F(A)\tens F(B)$.  
Following \cite[Section 3]{SS3}, we say that $(F,G)$ is a {\bf weak
monoidal Quillen equivalence} if $G$ is lax monoidal and two extra
conditions hold.  First, for some cofibrant replacement $A\ra
\unit_\cC$, the induced map $F(A) \ra F(\unit_\cC)\ra \unit_\cD$ is a
weak equivalence.  Second, for any two cofibrant objects $A,B\in \cC$
the map $F(A\tens B) \ra F(A)\tens F(B)$ is a weak equivalence.


\section{Adjoint modules}\label{ap:adj.modules}
\label{se:adjmod}

In this section and the next we deal with the general situation of one
Quillen  pair enriched over another Quillen pair.  Let
$\cC$ and $\cD$ be symmetric monoidal model categories, let $\cM$ be a
$\cC$-model category, and let $\cN$ be a $\cD$-model category.  Let
\[ F\colon \cC \adjoint \cD \colon G \qquad\text{and}\qquad
L\colon \cM \adjoint \cN \colon R 
\]
be two Quillen pairs, where we assume that $G$ is lax monoidal (see
Section~\ref{se:monoidalfunc}).  As usual, we'll write $\mM_{\cC}(X,Y)$
and $\mN_\cD(X,Y)$ for the enriched morphism objects over $\cC$ and
$\cD$, respectively.
 
Finally, let $Y$ be a cofibrant-fibrant object in $\cN$.  Then
$\mN_\cD(Y,Y)$ is a monoid in $\cD$, and so $G(\mN_\cD(Y,Y))$ is a
monoid in $\cC$.  Alternatively, we may choose a cofibrant-replacement
$\cQ RY \trfib RY$ and consider the $\cC$-monoid $\mM_\cC(\cQ RY,\cQ
RY)$.  How can we compare these two monoids, and under what conditions
will they be weakly equivalent?

This question can be answered by requiring certain compatibility
conditions between $(L,R)$ and $(F,G)$.  The goal of the
present section is to write down these conditions; this culminates in
Definition~\ref{de:Qmodule}, where we define what it means for $(L,R)$
to be an \dfn{adjoint module} over $(F,G)$.  The next section uses
this to tackle the problem of comparing enrichments.

\medskip
 
\subsection{Compatibility structure}
 
Before we can develop the definition of an adjoint module
we need the following statement. 
For the moment we only assume that $(F,G)$ and $(L,R)$ are
adjunctions.  That is, we temporarily drop the assumptions that they
are Quillen pairs and that $G$ is lax monoidal.
 
\begin{prop}
\label{pr:gamma}
There is a canonical bijection between natural transformations of the
following four types:
\begin{enumerate}[(i)]
\item $G\mN_\cD(LX,Y) \ra \mM_\cC(X,RY)$
\item $L(X\tens c) \ra LX \tens Fc$
\item $RY\tens Gd \ra R(Y\tens d)$
\item $G\mN_\cD(X,Y) \ra \mM_\cC(RX,RY)$.
\end{enumerate}
\end{prop}
 
\begin{proof}
This is a routine exercise in adjunctions.  We will only do some
pieces of the argument and leave the rest to the reader.
 
Suppose given a natural transformation $G\mN_\cD(LX,Y) \ra
\mM_\cC(X,RY)$.  For any $c\in \cC$ one therefore has the composition
\begin{myequation}
\label{eq:gamma}
 \xymatrixcolsep{0.9pc}\xymatrix{
& \cC(c,G\mN_\cD(LX,Y)) \ar[d]^\iso\ar[r] & \cC(c,\mM_\cC(X,RY)) \ar[d]^\iso \\\cN(LX\tens Fc,Y) \ar[r]^-\iso
& \cD(Fc,\mN_\cD(LX,Y))
&\cM(X\tens c,RY) \ar[r]^-\iso & \cN(L(X\tens c),Y).
}
\end{myequation}
By the Yoneda Lemma this gives a map
$L(X\tens c) \ra LX\tens Fc$, and this is natural in both $X$ and $c$.
 
Likewise, suppose given a natural transformation $L(X\tens c) \ra
LX\tens Fc$.  Then for $Y\in \cN$ and $d\in \cD$ we obtain
\[ L(RY\tens Gd) \ra LRY \tens FGd \ra Y \tens d \]
where the  second map uses the units of the adjunctions.
Taking the adjoint of the composition gives $RY\tens Gd \ra R(Y\tens
d)$, as desired.
 
Finally, suppose again that we have a natural transformation
$G\mN_\cD(LX,Y) \ra \mM_\cC(X,RY)$.  For $X,Y\in \cN$ consider the
composite
\[ G\mN_\cD(X,Y) \ra G\mN_\cD(LRX,Y) \ra \mM_\cC(RX,RY) \]
where the first map is obtained by applying $G$ to $\mN_\cD(X,Y) \ra
\mN_\cD(LRX,Y)$ induced by the unit $LRX\ra X$.  The above composite
is our natural transformation of type (iv).
 
We have constructed maps $(i)\ra (ii)$, $(ii)\ra (iii)$, and $(i)\ra
(iv)$.  We leave it to the reader to construct maps in the other
directions and verify that one obtains inverse bijections.
\end{proof}

\begin{remark}
Suppose we are given a natural transformation $\gamma\colon
G\mN_\cD(LX,Y) \ra \mM_\cC(X,RY)$.  Using the bijections from the
above result, we obtain natural transformations of types (ii), (iii),
and (iv).  We will also call each of these $\gamma$, by abuse.
\end{remark}

The next proposition lists the key homotopical properties 
required for $(L,R)$ to be a Quillen adjoint module over 
$(F,G)$.  
 
\begin{prop}
\label{pr:weconds}
Assume that $(F,G)$ and $(L,R)$ are Quillen pairs and that $\gamma\colon
G\mN_\cD(LX,Y) \ra \mM_\cC(X,RY)$ is a natural transformation.
\begin{enumerate}[(a)]
\item
The following two conditions are equivalent:
\begin{itemize}
\item The map $\gamma\colon G\mN_\cD(LX,Y) \ra \mM_\cC(X,RY)$ is a weak
equivalence whenever $X$ is cofibrant and $Y$ is fibrant.
\item The map $\gamma\colon L(X\tens c) \ra LX \tens Fc$ is a weak
equivalence whenever $X$ and $c$ are both cofibrant.
\end{itemize}
\item If $(L,R)$ is a Quillen equivalence, the conditions in (a)
are also equivalent to:
\begin{itemize}
\item For any cofibrant replacement $\cQ RX\ra RX$, the composite map
\[ G\mN_\cD(X,Y) \llra{\gamma} \mM_\cC(RX,RY) \ra \mM_\cC(\cQ RX,RY) \]
is a weak equivalence whenever $X$ is cofibrant-fibrant and $Y$ is fibrant.
\end{itemize}
\item Assume that both $(L,R)$ and $(F,G)$ are Quillen equivalences.
Then the conditions in (a) and (b) are also equivalent to:
\begin{itemize}
\item
For any cofibrant replacements $\cQ RY \ra RY$ and $\cQ' Gd \ra Gd$
and any fibrant replacement $Y\tens d \ra \ff(Y\tens d)$, the
composite
\[ \cQ RY \tens \cQ'Gd \ra RY\tens Gd \llra{\gamma} R(Y\tens d) \ra
R\ff(Y\tens d)\]
is a weak equivalence whenever $Y$ and $d$ are cofibrant and fibrant.
\end{itemize}
\end{enumerate}
\end{prop}
 
\begin{proof}
This is routine and basically follows from the adjunctions in 
Proposition~\ref{pr:gamma} with the following two additions.  For the 
equivalence in part (a),
consider the maps from~\ref{eq:gamma} in the respective homotopy categories.
For the equivalence with (b), note that the composite in (b) agrees with the 
composite
\[ G\mN_\cD(X,Y) \ra G\mN_\cD(L\cQ RX,Y) \llra{\gamma} \mM_\cC(\cQ RX,RY).\]
\end{proof}
 
The above homotopical properties need to be supplemented by
categorical associativity and unital properties which are listed in
the next two propositions.  Then, after stating these categorical
properties, we finally state the definition of a Quillen adjoint
module.
 
\begin{prop}
\label{pr:module}
Assume $G$ is lax monoidal.  Note that this gives a lax comonoidal
structure on $F$, by adjointness.  Let $\gamma$ again denote a set of
four corresponding natural transformations of types (i)--(iv).
Then the conditions in (a) and (b) below are equivalent:
 
\begin{enumerate}[(a)]
\item
The diagrams
\[ \xymatrix{L((X\tens c)\tens c')
\ar[r]^\gamma\ar[d]^\iso & L(X\tens c) \tens Fc'
\ar[r]^{\gamma\tens 1} & (LX\tens Fc) \tens Fc' \ar[d]^\iso \\
L(X\tens(c\tens c')) \ar[r]^\gamma & LX\tens F(c\tens c') \ar[r] &
LX\tens (Fc\tens Fc')
}
\]
all commute, for any $X$, $c$, $c'$.
\item The diagrams
\[ \xymatrix{
RY \tens (Gd \tens Gd') \ar[r]\ar[d]^\iso &
 RY\tens G(d\tens d') \ar[r]^\gamma
& R(Y\tens (d\tens d')) \ar[d]^\iso \\
(RY\tens Gd)\tens Gd' \ar[r]^{\gamma\tens 1} & R(Y\tens d) \tens Gd'
 \ar[r]^\gamma &
R((Y\tens d)\tens d')
}
\]
all commute, for any $Y$, $d$, $d'$.
\end{enumerate}
If $G$ is lax symmetric monoidal, then the above (equivalent) conditions imply
the following one:
\begin{enumerate}[(a)]
\addtocounter{enumi}{2}
\item The diagrams
\[ \xymatrixcolsep{.8pc}\xymatrix{
G\mN_\cD(Y,Z)\tens G\mN_\cD(X,Y) \ar[r]\ar[d]_{\gamma\tens \gamma}
 & G\biggl ( \mN_\cD(Y,Z)\tens
\mN_\cD(X,Y) \biggr ) \ar[r] & G\mN_\cD(X,Z)\ar[d]^\gamma \\
\mM_\cC(RY,RZ) \tens \mM_\cC(RX,RY) \ar[rr] && \mM_\cC(RX,RZ)
}
\]
commute for any $X$, $Y$, and $Z$.
\end{enumerate}
\end{prop}
 
\begin{proof}
The equivalence of (a) and (b) is extremely tedious but routine; we
leave it to the reader.  For (c), note that by using the adjunction
$\cC(c,\mM_\cC(RX,RZ))\iso \mM_\cC(RX\tens c,RZ)$ the two ways of going
around the diagram correspond to two maps
\[ RX\tens [G\mN_\cD(Y,Z)\tens G\mN_\cD(X,Y)] \lra RZ.
\]
One of these is the composite
\[\xymatrixcolsep{1pc}\xymatrix{
& RX\tens [G\mN_\cD(Y,Z)\tens G\mN_\cD(X,Y)] \ar[r] &
RX\tens G[\mN_\cD(Y,Z)\tens \mN_\cD(X,Y)]  \ar[d] \\
RZ & R(X\tens \mN_\cD(X,Z))\ar[l]  & RX\tens G[\mN_\cD(X,Z)] \ar[l]^\gamma
}
\]
The other is the composite
\[ \xymatrixcolsep{1.2pc}\xymatrix{
RX\tens [G\mN_\cD(Y,Z)\tens G\mN_\cD(X,Y)] \ar[r]^\iso
& RX \tens [G\mN_\cD(X,Y)\tens G\mN_\cD(Y,Z)] \ar[d] \\
[ R(X\tens \mN_\cD(X,Y))] \tens G\mN_\cD(Y,Z)  \ar[d]
 & [ RX\tens G\mN_\cD(X,Y)] \tens G\mN_\cD(Y,Z) \ar[l]_{\gamma\tens 1} \\
RY\tens G\mN_\cD(Y,Z) \ar[r]^\gamma & R(Y\tens \mN_\cD(Y,Z)) \ar[r] & RZ.
}
\]
The commutativity isomorphism comes into the first stage of this
composite because of how the composition map $\mM_\cC(RY,RZ)\tens
\mM_\cC(RX,RY) \ra \mM_\cC(RX,RZ)$ relates to the evaluation maps
under adjunction---see \cite[Prop. A.3]{hend}, for instance.

It is now a tedious but routine exercise to prove that the
above two maps
\[ RX\tens [G\mN_\cD(Y,Z)\tens G\mN_\cD(X,Y)] \dbra RZ
\] 
are indeed the same.  One forms the adjoints and then writes down a
huge commutative diagram.  A very similar result (in fact, a special
case of the present one) is proven in \cite[A.9]{hend}.
\end{proof}
 
Note that if $G$ is lax monoidal then it comes with a prescribed map
$\unit_\cC \ra G(\unit_\cD)$; adjointing gives $F(\unit_\cC)\ra
\unit_\cD$.  The following result concerns compatibility between these
maps and $\gamma$:
 
\begin{prop}
\label{pr:unit}
Assume again that $G$ is lax monoidal, and let $\gamma$ denote a set
of four corresponding natural transformations of types (i)--(iv).
The following three conditions are equivalent:
\begin{enumerate}[(a)]
\item For any $X$, the following square commutes:
\[ \xymatrix{
LX \ar[r]^-\iso \ar[d]_\iso & LX\tens \unit_\cD  \\
L(X\tens \unit_\cC) \ar[r]^-\gamma & LX\tens F(\unit_\cC).\ar[u]}
\]
\item For any $Y$, the following square commutes:
\[ \xymatrix{
RY \ar[r]^\iso\ar[d]_\iso & RY\tens \unit_\cC \ar[d] \\
R(Y\tens \unit_\cD)  & RY\tens G(\unit_\cD).\ar[l]_\gamma}
\]
\item For any $Y$, the following square commutes:
\[ \xymatrix{ \unit_\cC \ar[r]\ar[d] & G(\unit_\cD) \ar[d] \\
\mM_\cC(RY,RY) & G\mN_\cD(Y,Y) \ar[l]_-\gamma}
\]
\end{enumerate}
\end{prop}
 
\begin{proof}
Left to the reader.
\end{proof}

Finally we have the main definition:

\begin{defn}
\label{de:Qmodule}
Assume given adjoint pairs $(F,G)$ and $(L,R)$ where $G$ is lax
monoidal.  We will say that $(L,R)$ is an \dfn{adjoint module} over
$(F,G)$ if there exists a natural transformation $\gamma\colon
L(X\tens c) \ra LX \tens Fc$ such that the conditions of
Propositions~\ref{pr:module}(a) and ~\ref{pr:unit}(a) are both
satisfied.
 
If in addition $(F,G)$ and $(L,R)$ are both Quillen pairs and the equivalent
conditions of Proposition~\ref{pr:weconds}(a) are satisfied we will say
that $(L,R)$ is a \dfn{Quillen adjoint module} over $(F,G)$.
\end{defn}

\subsection{Basic properties} 
Below we give three properties satisfied by Quillen adjoint modules.
Recall the notion of a $\cC$-Quillen
adjunction between $\cC$-model categories, as in \cite[A.7]{hend}.
This is a Quillen pair $L\colon \cM \adjoint \cN \colon R$ where $\cM$
and $\cN$ are $\cC$-model categories, together with natural
isomorphisms $L(X\tens c) \iso L(X)\tens c$ which reduce to the
canonical isomorphism for $c=\unit_\cC$ and which are compatible with
the associativity isomorphisms in $\cM$ and $\cN$.  See also
\cite[Def. 4.1.7]{hovey}.

\begin{prop}
Suppose  $\cM$ and $\cN$ are $\cC$-model categories and
$L\colon \cM \adjoint \cN\colon R$ is a $\cC$-Quillen adjunction.
Then $(L,R)$ is a Quillen adjoint module over the pair
$(\id_\cC,\id_\cC)$.
\end{prop}
 
\begin{proof}
Since $(L,R)$ is a $\cC$-adjunction, there are
natural isomorphisms $LX \tens c \to L(X\tens c)$ which satisfy the
associativity and unital properties listed in
Propositions~\ref{pr:module}(a) and~\ref{pr:unit}(a).  This also
fulfills the second condition listed in
Proposition~\ref{pr:weconds}(a)
\end{proof}
 
\begin{prop}
Let $F\colon\cC\adjoint \cD\colon G$ be a Quillen pair between
symmetric monoidal model categories, where $G$ is lax monoidal.  Let $F'\colon
\cD \adjoint \cE \colon G'$ be another such pair.  Let $L\colon \cM
\adjoint \cN \colon R$ and $L'\colon \cN\adjoint \cP \colon R'$ be
Quillen pairs such that $(L,R)$ is a Quillen adjoint module over
$(F,G)$ and $(L',R')$ is a Quillen adjoint module over $(F',G')$.
Then $(L'L,RR')$ is a Quillen adjoint module over $(F'F,GG')$.
\end{prop}
 
\begin{proof}
For $X\in \cM$ and $c\in \cC$ we have natural maps
\[ L'L(X\tens c) \ra L'(LX\tens Fc) \ra L'(LX) \tens F'(Fc)
\]
using the adjoint module structure on $(L,R)$ over $(F,G)$ first, and
the module structure on $(L',R')$ over $(F',G')$ second.  One just has
to check the axioms to see that these maps make $(L'L,RR')$ a Quillen
adjoint module over $(F'F,GG')$.  This is a routine exercise in
categorical diagramming which we will leave to the reader.
\end{proof}
 
\begin{cor}
\label{co:adjointcompose}
Suppose $(L,R)$ is a Quillen adjoint module over $(F,G)$, and also
suppose that $\cP$ is a $\cC$-model category and $J\colon \cP \adjoint
\cM \colon K$ is a $\cC$-Quillen adjunction.  Then $(LJ,KR)$ is a
Quillen adjoint module over $(F,G)$.
\end{cor}
 
\begin{proof}
This is an immediate consequence of the above two propositions.
\end{proof}


\section{Applications of adjoint modules}\label{ap:ap.adj} 
\label{se:adjmodapp}

Recall from the last section that $\cC$ and $\cD$ are symmetric
monoidal model categories, $\cM$ is a $\cC$-model category, and $\cN$
is a $\cD$-model category.  We have Quillen pairs
\[ F\colon \cC \adjoint \cD \colon G \qquad\text{and}\qquad
L\colon \cM \adjoint \cN \colon R 
\]
in which $G$ is lax monoidal, and we assume that $(L,R)$ is a Quillen
adjoint module over $(F,G)$ as defined in Definition~\ref{de:Qmodule}.

Recall the notion of model enrichment from
Section~\ref{se:modelenrich}.  The assignment $X,Y \mapsto
\mN_\cD(X,Y)$ is a $\cD$-model enrichment of $\cN$, as in
\cite[Example 3.2]{hend}.  The induced assignment $X,Y\mapsto
G\mN_\cD(X,Y)$ is a $\cC$-model enrichment of $\cN$, by
Proposition~\ref{prop-mod-en} below.  Alternatively, if $\cQ W\trfib W$
is a cofibrant-replacement functor for $\cM$ and $W\trcof \ff W$ is a
fibrant-replacement functor for $\cN$, then one obtains another
$\cC$-model enrichment of $\cN$ via $X,Y \mapsto \mM_\cC(\cQ R\ff X,
\cQ R\ff Y)$.  (This is precisely the enrichment $L_*[\mM_\cC]$, as
defined in \cite[Section 3.4]{hend}.)  If $R$ preserves all weak
equivalences, the simpler assignment $X,Y \mapsto \mM_\cC(\cQ RX, \cQ RY)$
is also a $\cC$-model enrichment.

\begin{thm}
\label{th:main}
Assume the pair $(L,R)$ is a Quillen adjoint module over $(F,G)$.
Also assume that $G$ is lax symmetric monoidal and that $(L,R)$ is a
Quillen equivalence.  Then the two $\cC$-model enrichments on $\cN$
given by $X,Y \mapsto G\mN_\cD(X,Y)$ and $X,Y \mapsto \mM_\cC(\cQ R\ff X,
\cQ R\ff Y)$ are quasi-equivalent.  That is to say, $L_*\mM_\cC \he
G\mN_\cD$.  

If $R$ preserves all weak equivalences, then the above enrichments are
also quasi-equivalent to
$X,Y \mapsto \mM_\cC(\cQ R X, \cQ R Y)$.
\end{thm}
 
The above theorem  compares enrichments which have been transferred
over the right adjoints.  We would like to consider transfers over
left adjoints as well.  The situation is not completely dualizable,
though.  This is because there are no general conditions which ensure
$\mM_\cC(X,Y)$ is cofibrant, and so $F\mM_\cC(X,Y)$ will usually not
have the correct homotopy type.  

We do have the following corollary, however:
 
\begin{cor}
\label{co:GS}
Under the assumptions of the theorem, the two $\cC$-model enrichments
on $\cM$ given by $X,Y \mapsto G\mN_\cD(\ff L \cQ X,\ff L \cQ Y)$ and
$X,Y \mapsto \mM_\cC(X, Y)$ are quasi-equivalent.  That is,
$L^*[G\mN_{\cD}]$ is quasi-equivalent to $\mM_{\cC}$.  
\end{cor}
 
The quasi-equivalences in Theorem~\ref{th:main} and
Corollary~\ref{co:GS} are used in a key argument in~\cite{GS} to
translate a construction in $H\Q$-algebras into rational dgas.  The
following immediate corollary of the above theorem is what we will
mainly need in the present paper.
 
\begin{cor}
\label{co:main}
Assume that $\cC$ is combinatorial, satisfies the monoid axiom, and
that $\unit_\cC$ is cofibrant.
Under the assumptions of the theorem,
let $X\in \cN$ be a cofibrant-fibrant object.  Let $A\in \cM$ be any
cofibrant-fibrant object which is weakly equivalent to $RX$.
Then the $\cC$-monoids
$G\mN_\cD(X,X)$ and $\mM_\cC(A,A)$ are weakly equivalent.
\end{cor}
 
The extra assumptions on $\cC$ are necessary in order to apply a
certain proposition from \cite{hend}, saying that quasi-equivalent
enrichments give weakly equivalent endomorphism monoids.  
 
\subsection{Proofs of the above results}

\begin{prop}\label{prop-mod-en}
The assignment $n,n' \mapsto G\mN_\cD(n,n') $ is a $\cC$-model
enrichment on $\cN$.
\end{prop}
 
\begin{proof}
One uses the monoidal structure on $G$ to produce the associative
and unital composition maps.  Since $G$ preserves equivalences between
all fibrant objects and $\mN_{\cD}(n,n')$ is fibrant
if $n$ is cofibrant and $n'$ is fibrant, we see that
$G\mN_{\cD}(a',x) \to G\mN_{\cD}(a,x)$ and
$G\mN_{\cD}(a,x) \to G\mN_{\cD}(a,x')$ are weak equivalences whenever
$a \to a'$ is a weak equivalence between cofibrant objects and $x \to x'$
is a weak equivalence between fibrant objects.
\end{proof}

\begin{proof}[Proof of Theorem~\ref{th:main}]
For $X,Y\in \cN$ define
 $\sigma(X,Y)=G\mN_\cD(\ff X,\ff Y)$ and $\tau(X,Y)=\mM_\cC(QR\ff
X,QR\ff Y)$.  These are both $\cC$-model enrichments on $\cN$, and the
former is quasi-equivalent to $X,Y\mapsto G\mN_\cD(X,Y)$ by
\cite[Prop. 3.9]{hend}.
 
Define $W(X,Y)=\mM_\cC(\cQ R\ff X,R\ff Y)$.  This is a $\sigma-\tau$
bimodule via the maps
\[\xymatrixcolsep{1.4pc}\xymatrix{
G\mN_\cD(\ff Y,\ff Z) \tens
\mM_\cC(\cQ R\ff X, R\ff Y) \ar[r]^-{\gamma \tens 1}
 & \mM_\cC(R\ff Y,R\ff Z)\tens
\mM_\cC(\cQ R\ff X,R\ff Y) \ar[d] \\
 & \mM_\cC(\cQ R\ff X,R\ff Z) }
\]
and
\[ \mM_\cC(\cQ R\ff Y,R\ff Z) \tens \mM_\cC(\cQ R\ff X,\cQ R\ff Y) \lra
\mM_\cC(\cQ R\ff X,R\ff Z).
\]
Some routine but tedious checking is required to see that this indeed
satisfies the bimodule axioms of \cite[Section 2.2]{hend}.  This uses
the conditions from Proposition~\ref{pr:module}(c) and
Proposition~\ref{pr:unit}(a).
 
The canonical maps $\cQ R\ff X\ra R\ff X$ give maps $\unit_{\cC} \ra W(X,X)$
making $W$ into a pointed bimodule, and one checks using the
condition from Proposition~\ref{pr:weconds}(b) that this is a
quasi-equivalence.  This last
step uses our assumption that $(L,R)$ is a Quillen equivalence.
 
If $R$ preserves all weak equivalences, then the above proof works
even if every appearance of the functor $\ff$ is removed.
\end{proof}
 
\begin{proof}[Proof of Corollary~\ref{co:GS}]
The result \cite[3.14(d)]{hend} shows that since $L$ is a Quillen
equivalence the maps $L^*$ and $L_*$ are inverse bijections.  Since we
have already proven $L_*\mM_{\cC}\he G\mN_{\cD}$, we must have
$L^*[G\mN_{\cD}] \he \mM_{\cC}$.
\end{proof}
 
\begin{proof}[Proof of Corollary~\ref{co:main}]
Using the above theorem together with \cite[Cor. 3.6]{hend} (which
requires our assumptions on $\cC$) we find that if
$X\in \cN$ is cofibrant-fibrant then the $\cC$-monoids $G\mN_\cD(X,X)$
and $\mM_\cC(\cQ R\ff X, \cQ R\ff X)$ are weakly equivalent.  However, note
that one has a weak equivalence $A \we \cQ R\ff X$.  By applying
\cite[Cor. 3.7]{hend} (in the case where $\cI$ is the category with
one object and an identity map) one finds that the $\cC$-monoids
$\mM_\cC(\cQ R\ff X,\cQ R\ff X)$ and $\mM_\cC(A,A)$ are weakly equivalent.
\end{proof}
 
 
\subsection{Applications to module categories}
\label{se:app-module-cat}
We'll now apply the above results to the homotopy theory of $\cC
I$-categories.  Readers may want to review Appendix~\ref{se:CI-cat}
before proceeding further.

Let $\cC$ and $\cD$ be cofibrantly-generated symmetric monoidal model
categories satisfying the monoid axiom, and assume that $\unit_\cC$
and $\unit_\cD$ are cofibrant.  Let $F\colon \cC \adjoint \cD \colon
G$ be a Quillen pair where $G$ is lax monoidal.  Let $I$ be a set and
consider the notion of $\cC I$-category (a category enriched over
$\cC$ with object set $I$) from Appendix~\ref{se:CI-cat}.  Note that when $I$
consists of one object then a $\cC I$-category is just a monoid in
$\cC$.

Let $\cR$ be a $\cD I$-category, and consider the category $\Mod \cR$
of right $\cR$-modules.  By~\cite[6.1]{SS3} 
the category $\Mod \cR$ has a model structure in which the weak
equivalences and fibrations are obtained by forgetting objectwise to
$\cD$.  This is a $\cD$-model category in a natural way.  The SM7 (or
pushout product) condition follows from $\cD$ using~\cite[3.5]{SS1}
since the $\cD$ action is pointwise, and the unit condition follows
from our assumption that $\unit_\cD$ is cofibrant (since this implies
that the cofibrant $\cR$-modules are objectwise cofibrant).
   
Since $G$ is lax monoidal, $G\cR$ is a $\cC I$-category and we may
consider the corresponding module category $\Mod G\cR$.  This is a
$\cC$-model category.  If $M$ is an
$\cR$-module then $GM$ becomes a $G\cR$-module in a natural way, and there
is an adjoint pair $F^{\cR} \colon \Mod G\cR \adjoint \Mod \cR \colon
G$ by Proposition~\ref{pr:CImod}(a).
The functors $(F^\cR,G)$ are a Quillen pair since $G$ preserves the
objectwise fibrations and trivial fibrations.
 
We are now in the position of having two Quillen pairs $F\colon
\cC\adjoint \cD \colon G$ and $F^{\cR}\colon \Mod G\cR \adjoint \Mod\cR
\colon G$.  The categories $\Mod G\cR$ and $\Mod \cR$ are $\cC$- and
$\cD$-model categories, respectively.

\begin{prop}
\label{pr:mod1}
Under the above assumptions on $\cC$, $\cD$, and $G$ one has:
\begin{enumerate}[(a)]
\item $(F^{\cR},G)$ is an adjoint module over $(F,G)$.
\item If $F$ is strong monoidal, then $(F^{\cR},G)$ is a Quillen
adjoint module over $(F,G)$.
\item Assume that $\cC$ is a stable model category whose homotopy
category is generated by $\unit_\cC$.  Assume as well that
$F(\unit_\cC)\to \unit_\cD$ is a weak equivalence. 
Then $(F^{\cR},G)$ is a Quillen adjoint module
over $(F,G)$.
\end{enumerate}
\end{prop}

\begin{proof}
In terms of the notation of Section~\ref{se:adjmod} we have $L=F^{\cR}$
and $R=G$.  A natural transformation $\gamma$ of the type in
Proposition~\ref{pr:gamma}(iii) is therefore obtained using the lax
monoidal structure on $G$.  This automatically satisfies the axioms of
Proposition~\ref{pr:module}(b) and Proposition~\ref{pr:unit}(b), so
that we have an adjoint module over $(F,G)$.  This proves (a).
 
To prove (b) we show that $L(X \tens c) \to LX \tens Fc$ is an isomorphism,
and hence a weak equivalence.  Here $L = F^{\cR} = F (-) \tens_{FG\cR} \cR$
since $F$ is strong monoidal; see the discussion above~\cite[3.11]{SS3}.  
  It is then easy to verify that $F^{\cR}(X \tens c)=
F(X \tens c) \tens_{FG\cR} \cR \iso (FX \tens Fc) \tens_{FG\cR} \cR
\iso F^{\cR}(X) \tens Fc$.
 
To prove (c), we will verify that $G\mN_\cD(LX,Y)\llra{\gamma}
\mM_\cC(X,RY)$ is a weak equivalence whenever $X$ is cofibrant and $Y$
is fibrant.  Using our assumption about $\unit_\cC$ generating
$\ho(\cC)$, it suffices to show that
\[ [\unit_\cC,G\mN_\cD(LX,Y)]_* \ra [\unit_\cC,\mM_\cC(X,RY)]_* \]
is an isomorphism of graded groups, where $[\blank,\blank]_*$ denotes
the graded group of maps in a triangulated category.

By adjointness, the problem reduces to showing that the map 
\[ [LX\tens F( 
\unit_\cC),Y]_* \ra [L(X\tens \unit_\cC),Y]_*
\] 
is an isomorphism---or
in other words, that $LX\tens F(\unit_\cC) \ra L(X\tens \unit_\cC)$
is a weak equivalence.  But this follows easily from our assumption
that $F(\unit_\cC)\ra\unit_\cD$ is a weak equivalence.
\end{proof}
 
Now assume that $\cO$ is a cofibrant $\cC I$-category.
By Proposition~\ref{pr:CIcat} there is an adjunction
$F^{\cD I}\colon \cC I-\Cat \adjoint \cD I-\Cat \colon G$, so that we
get a $\cD I$-category $F^{\cD I} \cO$.
By Proposition~\ref{pr:CImod}(b) there is a
Quillen pair
\[ F_{\cO} \colon \Mod \cO \adjoint \Mod F^{\cD I} \cO \colon G_\cO.
\]
 
\begin{prop}
\label{pr:mod2}
In the above setting one has:
\begin{enumerate}[(a)]
\item $(F_{\cO},G_{\cO})$ is an adjoint module over $(F,G)$.
\item If $F$ is strong monoidal,
then $(F_{\cO},G_{\cO})$ is a Quillen adjoint module over $(F,G)$.
\item Assume that $\cC$ is a stable model category whose homotopy
category is generated by $\unit_\cC$, and that
$F(\unit_\cC)\to \unit_\cD$ is a weak equivalence.
Then $(F_{\cO},G_{\cO})$ is a Quillen adjoint module over $(F,G)$.
\end{enumerate}
\end{prop}
 
\begin{proof}
Write $\cR=F^{\cD I}\cO$.
The adjunction $(F_{\cO},G_{\cO})$ is the composite of the two adjunctions
\[ \xymatrix{
\Mod \cO  \ar@<0.5ex>[r]^-{\beta_*} & \Mod G\cR
\ar@<0.5ex>[r]^-{F^{\cR}}
\ar@<0.5ex>[l]^-{\beta^*}
& \Mod \cR \ar@<0.5ex>[l]^-G
}
\]
where $\beta\colon \cO \ra G\cR=GF^{\cD I}\cO$ is the unit of the
adjunction $(F^{\cD I}, G)$.
  
But $(\beta_*,\beta^*)$ is a $\cC$-Quillen adjunction
and by Proposition~\ref{pr:mod1}, under either set of conditions,
we know $(F^\cR,G)$ is a Quillen adjoint
module over $(F,G)$.  The result now follows immediately from
Corollary~\ref{co:adjointcompose}.
\end{proof}

\begin{cor}\label{cor:mod2}
In addition to our previous assumptions, assume that $G$ is 
lax symmetric monoidal
and $\cO$ is a cofibrant $\cC I$-category.  Suppose also that
($F_\cO$, $G_\cO$) is a Quillen equivalence and the hypotheses in
either part (b) or (c) hold from Proposition~\ref{pr:mod2}.  Let $X\in
\Mod (F^{\cD I}\cO)$ be a cofibrant-fibrant object and let $A\in
\Mod\cO$ be any module weakly equivalent to $G_{\cO}X$.  Then the
$\cC$-monoids
\[ \underline{\Mod\cO}_\cC(A,A) \qquad \text{and} \qquad
G\bigl [\underline{\Mod (F^{\cD I}\cO)}_{\cD}(X,X) \bigr ]
\]
are weakly equivalent.
\end{cor}

\begin{proof}
This follows from the above proposition and Corollary~\ref{co:main}.
\end{proof}
 
\begin{example}\label{ex1}
The adjoint pair $L\colon \Sp^\Sigma(\ch_+)\adjoint
\Sp^\Sigma(s\Ab)\colon \nut$ from \cite[4.3]{S} forms one example for
$(F,G)$.  The result \cite[3.4]{S}
shows that $\Sp^\Sigma(\ch_+)$ and $\Sp^\Sigma(s\Ab)$ are cofibrantly
generated symmetric monoidal model categories which satisfy the monoid
axiom.  The conditions in Proposition~\ref{pr:mod1}(c)
or~\ref{pr:mod2}(c) are verified in the last paragraph of the proof
of~\cite[4.3]{S}. 
Note, though, that $L$ is not strong monoidal.
This failure is due to the fact that the adjunction $N\colon s\Ab
\adjoint \ch_+ \colon \Gamma$ is not monoidal~\cite[2.14]{SS3}.
 
Corollary~\ref{cor:mod2} holds for $(L,\nut)$ in place of $(F,G)$
because $N$ is lax symmetric monoidal,
so its prolongation and $\nut$ are also lax symmetric monoidal.
The fact that $(L_\cO, (\nut)_\cO)$ is a Quillen equivalence follows
from~\cite[3.4, 4.3]{S} 
and~\cite[6.5(1)]{SS3}. 
See also Proposition~\ref{pr:CImod}(c). 
\end{example}

\section{Transporting enrichments}\label{se:transfer}

In this section we prove a technical result about transporting
enrichments.  This will be needed later, in the proof of
Proposition~\ref{pr:ch->sAb}.  The basic idea is as follows.  Suppose
$\cM$ is a $\cC$-model category, where $\cC$ is a certain symmetric
monoidal model category.  Assume also that $\cD$ is another symmteric
monoidal model category, and that one has a Quillen equivalence $\cC
\adjoint \cD$ which is compatible with the monoidal structure.  Then one
might hope to find a $\cD$-model category $\cN$ which is Quillen
equivalent to $\cM$, and where the Quillen equivalence aligns the
$\cC$- and $\cD$-structures.  In this section we prove one theorem
along these lines, assuming several hypotheses on the given data.

\medskip

We begin with the following two definitions:

\begin{defn}\label{def:CompGen}
Let $\cT$ be a triangulated category with infinite
coproducts.
\begin{enumerate}[(a)]
\item
An object $P\in \cT$ is called \dfn{compact} if
$\oplus_\alpha \cT(P,X_\alpha) \ra \cT(P,\oplus_\alpha
X_\alpha)$ is an isomorphism for every set of objects $\{X_\alpha\}$;
\item A set of objects $S\subseteq \cT$ is a \dfn{generating set} if
the only full, triangulated subcategory of $\cT$ which contains $S$
and is closed under arbitrary coproducts is $\cT$ itself.
If $S$ is a singleton set $\{P\}$ we say that $P$ is a
\dfn{generator}.
\end{enumerate}
\end{defn}

When $\cM$ is a stable model category we will call an object compact
if it is compact in $\Ho(\cM)$, and similarly for the notion of
generating set.
Most stable model categories of interest have a generating set of
compact objects.  For example, Hovey shows in~\cite[7.4.4]{hovey}
that this is true for 
any finitely-generated, stable model category.

Let $\cC$ and $\cD$ be symmetric monoidal, stable model categories.
Let $\cM$ be a pointed $\cC$-model category (so that $\cM$ is also stable).
We make the following assumptions:
\begin{enumerate}[(a)]
\item $\cC$ and $\cD$ are combinatorial model categories satisfying
the monoid axiom, and their units are cofibrant.
\item There is a weak monoidal Quillen equivalence $F\colon \cC
\adjoint \cD \colon G$, where $G$ is lax symmetric monoidal.
\item $\cC$ satisfies axioms (QI1-2) from Appendix~\ref{se:CI-cat}.
\item
$\cC$ is a stable model category whose homotopy
category is generated by $\unit_\cC$, and
$F(\unit_\cC)\to \unit_\cD$ is a weak equivalence.
\item $\cM$ has a generating set of compact objects.
\end{enumerate}

If $\cN$ is a $\cD$-model category, let $G\mN_\cD$
denote the assignment $X,Y\mapsto G\mN_{\cD}(X,Y)$.  By
Proposition~\ref{prop-mod-en} this
is a $\cC$-model enrichment of $\cN$.

\begin{prop}
\label{pr:technical}
Under the above conditions there exists a combinatorial $\cD$-model
category $\cN$ and a zig-zag of Quillen equivalences
\[ \cM \llla{L_1} \cM_1 \llra{L_2} \cN\]
such that the model enrichment $\mM_{\cC}$ is quasi-equivalent to
 $(L_1)_*(L_2)^*[G\mN_{\cD}]$.  

If, in addition, $\cC$ and $\cD$ are additive model categories (see
the following section for the definition)  then
$\cM_1$ and $\cN$ may also be chosen to be additive.
\end{prop}

By \cite[3.6]{hend} 
this yields the following immediate corollary:

\begin{cor}\label{cor:technical}
If $Y\in \cN$ and $X\in \cM$ are cofibrant-fibrant objects and $Y$ is
the image of $X$ under the derived functors of the above Quillen
equivalence $\cM \he \cN$, then the $\cC$-monoids $\mM_{\cC}(X,X)$ and
$G\mN_{\cD}(Y,Y)$ are weakly equivalent.
\end{cor}

\begin{proof}[Proof of Proposition~\ref{pr:technical}]
Constructing the model category $\cN$ will require several steps, and
we will start by just giving a sketch---then we will come back and
provide detailed justifications afterwards.  

Let $I$ denote a set of cofibrant-fibrant, compact objects which
generate $\cM$.  Let $\cO$ be the $\cC I$-category \cite[6.2]{Bor}
defined by $\cO(i,j)=\mM_{\cC}(i,j)$.  Then there is a $\cC$-Quillen
equivalence
\begin{myequation}
\label{eq:step1}
T\colon\Mod\cO \adjoint \cM\colon S 
\end{myequation}
where $\Mod\cO$ is the model category of right $\cO$-modules
(see Proposition~\ref{pr:Omod}).

Let $g\colon\ccO \ra \cO$ be a cofibrant-replacement for $\cO$ in the model
category of $\cC I$-categories (Proposition \ref{pr:CIcat}(a)).  Then
tensoring and restricting give the left and right adjoints of a
$\cC$-Quillen equivalence 
\begin{myequation}
\label{eq:step1.5}
g_*\colon \Mod \ccO \adjoint \Mod \cO \colon g^*
\end{myequation}
(see Proposition \ref{pr:Omod}(c)).

Next we use the functor $L^{\cD I}$ from
Proposition~\ref{pr:CIcat}(b).  This gives us a $\cD I$-category
$L^{\cD I}\ccO$ and a Quillen equivalence
\begin{myequation}
\label{eq:step2}
L_{\ccO}\colon \Mod \ccO \adjoint \Mod(L^{\cD I}\ccO) \colon \nut.
\end{myequation}
Let $\cN = \Mod(L^{\cD I}\ccO)$.  This is a $\cD$-model category, and
we have established a zig-zag of Quillen equivalences
\[ \cM \bwe \Mod \cO \bwe \Mod \ccO \we
\Mod(L^{\cD I}\ccO)=\cN.
\]
We set $\cM_1=\Mod\ccO$.  Note that if $\cC$ and $\cD$ are additive
model categories then by Corollary~\ref{co:Cmodel=additive} so are
$\cM_1$ and $\cN$ (since $\cM_1$ is a $\cC$-model category and $\cN$
is a $\cD$-model category).
 
\vspace{0.2in}

Now we fill in the details of the above sketch. 
The category of right
$\cO$-modules $\Mod\cO$ is defined in \cite[Section 6]{SS3}, and the
model structure on $\Mod\cO$ is provided in \cite[6.1(1)]{SS3}. 
See Appendix~\ref{sec-ci-cat} for a review.
To justify the Quillen equivalence in~(\ref{eq:step1}),  
define $S\colon \cM \ra \Mod\cO$ by letting $S(Z)$ be the functor
$i\mapsto \ehom_{\cC}(i,Z)$.  This obviously comes equipped with a
structure of right $\cO$-module.  The construction of the left adjoint
can be copied almost verbatim from \cite[3.9.3(i)]{SS2}, 
 which
handled the case where $\cC$ was $\Sp^\Sigma$.  The right adjoint
obviously preserves fibrations and trivial fibrations, so we have a
Quillen pair.  It is readily seen to be a $\cC$-Quillen pair.

Finally, that this is a Quillen equivalence follows just as in
\cite[3.9.3(ii)]{SS2}; 
this uses that $I$ was a generating set
of compact objects.  The proof can be summarized quickly as follows.
First, the compactness of the objects in $I$ shows that the derived
functor of $S$ preserves all coproducts; this is trivially true for
the derived functor of $T$ because it is a left adjoint.  One has
canonical generators $\Fr_i\in \Mod\cO$ for each $i\in I$, and
adjointness shows that $T(\Fr_i)\iso i$.  Likewise, $S(i)\iso \Fr_i$.
Using that the derived functors of $S$ and $T$ preserve coproducts and
triangles, one now deduces that the respective composites are naturally
isomorphic to the identities.
This completes step (\ref{eq:step1}) above.

We now turn to (\ref{eq:step1.5}).  The map of $\cC I$-categories
$g\colon \ccO \ra \cO$ gives a Quillen map $\Mod\ccO \ra \Mod\cO$ by
Proposition~\ref{pr:Omod}(b).  We will know this is a Quillen
equivalence by Proposition~\ref{pr:Omod}(c) as long as we know that
$\cC$ satisfies the axioms (QI1-2) of Appendix~\ref{sec-ci-cat}. 

The Quillen equivalence of (\ref{eq:step2}) is a direct application of
Proposition~\ref{pr:CImod}(c).

\medskip

At this point we have constructed the zig-zag $\cM \llla{L_1} \cM_1
\llra{L_2} \cN$.  We must verify that $(L_1)_*(L_2)^*[G\mN_{\cD}]$ is
quasi-equivalent to $\mM_{\cC}$.

It follows from Proposition~\ref{pr:mod2} and Theorem~\ref{th:main} that
$(L_2)^*[G\mN_{\cD}]$ is quasi-equivalent to $(\mM_1)_{\cC}$.  This is
where the theory of adjoint modules was needed.  Since
$L_1$ is a $\cC$-Quillen equivalence, it follows from
\cite[3.14(e)]{hend} 
that $(L_1)_*[(\mM_1)_{\cC}]$ is quasi-equivalent to $\mM_{\cC}$.  So
these two statements give exactly what we want.  
\end{proof}


\newpage

\section{Additive model categories}
\label{se:additivedef}

Now the second half of the paper begins.  We change direction and
start to pursue our main results on additive enrichments.  In the
present section we define the notion of an {\it additive\/} model
category, and prove some basic results for recognizing them.

\medskip

A category is \dfn{preadditive} if its hom-sets have natural
structures of abelian groups for which the composition pairing is
biadditive.  A category is \dfn{additive} if it is preadditive and
it has finite coproducts.  This forces the existence of an
initial object (the empty coproduct), which will necessarily be a zero
object.  See \cite[Section VIII.2]{MacL}.  A functor $F\colon \cC \ra
\cD$ between additive categories is an \dfn{additive functor} if
$F(f+g)=F(f)+F(g)$ for any two maps $f,g\colon X \ra Y$.

Now let $\cM$ be a model category whose underlying category is
additive.  Write $\cM_{cof}$ for the full subcategory of cofibrant
objects, and $c\cM$ for the category of cosimplicial objects in $\cM$.
Recall from \cite[Section 15.3]{H} that $c\cM$ has a Reedy model
category structure.  Also recall that a \dfn{cosimplicial resolution}
is a Reedy cofibrant object of $c\cM$ in which every coface and
codegeneracy map is a weak equivalence.

\begin{defn}
Let $\cI$ be a small, additive subcategory of $\cM_{cof}$.  
By an \dfn{additive cosimplicial resolution} on $\cI$ we mean an
additive functor $\Gamma\colon \cI \ra c\cM$ whose image lies in the
subcategory of cosimplicial resolutions, together with a natural
weak equivalence $\Gamma(X)^0\we X$.  
\end{defn}

By \cite[16.1.9]{H}, any small subcategory $\cI\subseteq \cM_{cof}$
has a cosimplicial resolution; however, the existence of an 
{\it additive\/} cosimplicial resolution is not at all clear.  

If $\Gamma$ and $\Gamma'$ are two additive cosimplicial resolutions on
$\cI$, then define a map $\Gamma\ra \Gamma'$ to be a natural transformation
of functors which gives commutative triangles
\[\xymatrix{
\Gamma(X)^0 \ar[d]\ar[r] & X \\
\Gamma'(X)^0 \ar[ur],
}
\]
for all $X\in \cI$.
The map is called a weak equivalence if all the maps $\Gamma(X) \ra
\Gamma'(X)$ are weak equivalences.

\begin{defn}
A model category $\cM$ is \dfn{additive} if its underlying category is
additive and if for every small, full subcategory $\cI$ of $\cM_{cof}$ the
following two statements are satisfied:
\begin{enumerate}[(a)]
\item $\cI$ has an additive cosimplicial resolution;
\item The category of additive cosimplicial resolutions on $\cI$,
where maps are natural weak equivalences, is connected (i.e., any two objects
are connected by a zig-zag).
\end{enumerate}
\end{defn}

\begin{remark}
One might argue that the adjective `connected' in the above definition
should be replaced with `contractible'.  This is a legitimate concern.
We have merely chosen the weakest definition which will support the
results in Section~\ref{se:univ}.
\end{remark}

\begin{prop}
\label{pr:additive-recognize}
Let $\cM$ be a model category whose underlying category is additive.
Suppose that there is a functor $F\colon \cM_{cof} \ra c\cM$
together with a natural isomorphism $F^0(X)\iso X$.  Assume that each
$F(X)$ is a cosimplicial resolution, that $F$ preserves colimits, and
that if $X\cof Y$ is a cofibration then $F(X) \ra F(Y)$ is a Reedy
cofibration.  Then $\cM$ is an additive model category.
\end{prop}

Note that the functor $F$ will automatically be additive; since it
preserves colimits, it preserves direct sums.

\begin{proof}
The existence of additive cosimplicial resolutions is provided by $F$.
So we must only prove that any two such resolutions can be connected
by a zig-zag.

If $\Gamma\in c\cM$ is any cosimplicial object, applying $F$ to
$\Gamma$ yields a bi-cosimplicial object $F\Gamma$ given by
$[m],[n]\mapsto F^m\Gamma^{n}$.  Let $\widetilde{\Gamma}\in c\cM$ denote
the diagonal of this bi-cosimplicial object, and note that there is a
natural map $\widetilde{\Gamma} \ra \Gamma$.  We claim that if
$\Gamma$ is a cosimplicial resolution then so is $\widetilde{\Gamma}$.

Suppose that $\Gamma\in c\cM$ is a cosimplicial resolution of some
object $X$.  Then every latching map $L^n\Gamma \ra \Gamma^n$ is a
cofibration (see \cite[15.3]{H} for a discussion of latching maps).
From the bi-cosimplicial object $F\Gamma$, we get a `vertical' 
latching map in
$c\cM$ of the form $L^{*,n}[F\Gamma]\ra F(\Gamma^n)$.  Here the domain
is the cosimplicial object which in level $m$ is the $n$th latching
object of $[F\Gamma]^{m,*}$.
Since the latching spaces are
formed as colimits, and $F$ preserves colimits, one has
$L^{*,n}[F\Gamma]\iso F(L^n\Gamma)$.  So we are looking at the map
$F(L^n\Gamma) \ra F(\Gamma^n)$.  But this is the result of applying
$F$ to a cofibration in $\cM$, so it is a Reedy cofibration.

So we are in the situation of Lemma~\ref{le:addlem} below, in which
every vertical latching map of $F\Gamma$ is a Reedy cofibration.  By
the lemma, this implies that the diagonal $\widetilde{\Gamma}$ is
Reedy cofibrant.  Since clearly every map in $F\Gamma$ is a weak
equivalence, it is therefore a cosimplicial resolution of $X$.

Now suppose that $\cI$ is a small, full subcategory of $\cM_{cof}$ and
$\Gamma_1,\Gamma_2\colon \cI \ra c\cM$ are two additive cosimplicial
resolutions.
For any $X\in \cI$ we have a canonical zig-zag $\Gamma_1(X) \we cX \bwe
\Gamma_2(X)$ where $cX$ denotes the constant cosimplicial object.
Consider the resulting diagram
\[ \xymatrix{
\widetilde{\Gamma}_1(X) \ar[r]^\sim \ar[d]^\sim & \widetilde{cX} \ar[d]^\sim &
\widetilde{\Gamma}_2(X) \ar[l]_\sim\ar[d]^\sim \\
\Gamma_1(X) \ar[r]^\sim & cX & \Gamma_2(X). \ar[l]_\sim
}
\]
The functors $\widetilde{\Gamma}_1,\widetilde{\Gamma}_2\colon \cI \ra
c\cM$ are additive cosimplicial resolutions on $\cI$.  So is the map
$\cI \ra c\cM$ given by $X\mapsto \widetilde{cX}=F(X)$.  Thus, the
outer rim of the
above diagram gives a zig-zag of weak equivalences connecting the
additive cosimplicial resolutions $\Gamma_1$ and $\Gamma_2$.  
\end{proof}

We need some notation for the following lemma.  Let $X^{*,*}$ be a
bi-cosimplicial object in a model category $\cM$.  Considering this as
an object of $c(c\cM)$, one obtains a `vertical' latching map $L^{*,n}X \ra
X^{*,n}$ in $c\cM$.  Here $L^{*,n}X$ denotes the cosimplicial object
sending $[m]$ to the $n$th latching object of $X^{m,*}$.  

\begin{lemma}
\label{le:addlem}
Let $\cM$ be any model category.  Suppose that $X^{*,*}$ is a
bi-cosimplicial object of $\cM$---that is, $X\in c(c\cM)$.  Assume
that every latching map $L^{*,n}X \ra X^{*,n}$ is a Reedy cofibration
in $c\cM$.  Then the diagonal cosimplicial object $[n]\mapsto X^{n,n}$
is Reedy cofibrant.
\end{lemma}

The proof of the above lemma is a little technical.  We defer it until
the end of the section.

\begin{cor}
\label{co:cor1}
Let $\cC$ and $\cM$ be model categories, where the underlying
category of $\cM$ is additive.  Suppose there is a bifunctor
$\tens\colon \cM \times \cC \ra \cM$ satisfying the pushout-product
axiom for cofibrations: if $i\colon A\cof B$ is a cofibration in $\cM$
and $j\colon X\cof Y$ is a cofibration in $\cC$, then $(A\tens Y)
\amalg_{(A\tens X)} (B\tens X) \ra B\tens Y$ is a cofibration which is
a weak equivalence if either $i$ of $j$ is.  Suppose also that
\begin{enumerate}[(i)]
\item
For
any $X\in \cC$ the functor $(\blank) \tens X$ preserves colimits;
\item 
For any $A\in \cM$ the functor $A\tens (\blank)$ preserves colimits;
\item There is a cofibrant object $\unit\in \cC$ and 
 natural isomorphisms $A\tens \unit \iso A$.
\end{enumerate}
Then $\cM$ is an additive model category.
\end{cor}

\begin{proof}
Let $\Gamma\in c\cC$ be a cosimplicial resolution of $\unit$ with
$\Gamma^0=\unit$.  For any cofibrant object $A\in \cM$, let $F(A)$ be
the cosimplicial object $[n]\mapsto A\tens \Gamma^n$.  The
pushout-product axiom, together with assumption (ii), shows that
$F(A)$ is a cosimplicial resolution of $A$.  Assumption (i) implies
that $F$ preserves colimits, and assumption (iii) says there are
natural isomorphisms $F(A)^0\iso A$.  Finally, it is an easy exercise
to use assumption (ii) and the pushout-product axiom to show that if
$A\ra B$ is a cofibration then $F(A)\ra F(B)$ is a Reedy cofibration.
The result now follows by applying
Proposition~\ref{pr:additive-recognize}.
\end{proof}

The above corollary lets one identify many examples of additive model
categories.  We only take note of the few obvious ones:

\begin{cor}
If $R$ is a ring, consider the model category $s(R-\dMod)$ where
fibrations and weak equivalences are determined by the forgetful
functor to $\sSet$.  This is an additive model category.  So is the
model category $\Ch(R)$ of unbounded chain complexes, where weak
equivalences are quasi-isomorphisms and fibrations are sujections.
\end{cor}

\begin{proof}
This results from two applications of the previous corollary.  For the
first statement we take $\cM=s(R-\dMod)$, $\cC=s(\Z-\dMod)$, and
$\tens$ to be the levelwise tensor product over $\Z$.  Here we are
using that if $M$ is an $R$-module and $A$ is a $\Z$-module then
$M\tens_\Z A$ has a natural $R$-module structure from the left.

For the second statement
we can take $\cM=\Ch(R)$, $\cC=\Ch_{\geq 0}(\Z)$, and $\tens$ the usual
tensor product of chain complexes over $\Z$.  (One could also take
$\cC=\Ch(\Z)$, but verifying the pushout-product axiom is a little
easier for bounded below complexes).
\end{proof}

If $R$ is a dga, then $R-\dMod$ has a model category structure where
weak equivalences are quasi-isomorphisms and fibrations are surjections.

\begin{cor}
If $R$ is a dga, then the model category $R-\dMod$ is additive.
\end{cor}

\begin{proof}
We again apply Corollary~\ref{co:cor1}, this time with $\cM=R-\dMod$
and $\cC=\Ch_{\geq 0}(\Z)$.  The $\tens$ functor is the tensor product
$M, C \mapsto M\tens_\Z C$ with the induced left $R$-module structure.
\end{proof}

We also note the following result:

\begin{cor}
\label{co:Cmodel=additive}
Let $\cC$ be a symmetric monoidal model category in which the unit is
cofibrant, and where the underlying category is additive.  Then $\cC$
is an additive model category.  Any $\cC$-model category is also
additive.
\end{cor}

\begin{proof}
The first statement follows immediately from Corollary~\ref{co:cor1}, as
the bifunctor $X,Z\mapsto X\tens Z$ preserves colimits in both
variables.

The second statement is also a direct application of
Corollary~\ref{co:cor1}, as soon as one notes that if $\cM$ is a
$\cC$-model category then the underlying category of $\cM$ is
additive.  This follows using the adjunctions $\cM(X,Y)\iso \cM(X\tens
\unit_\cC,Y)\iso \cC(\unit_\cC,\mM_{\cC}(X,Y))$, as there is a natural
abelian group structure on the latter set.  One checks that
composition is biadditive with respect to this structure.
\end{proof}

\subsection{Bisimplicial machinery}
The last thing we must do in this section is prove
Lemma~\ref{le:addlem}.  This requires some machinery which we briefly
recall.  

If $K\in \sSet$ and $A\in c\cM$, one may form the coend $A\tens K \in \cM$.
This is the coequalizer of the two arrows
\begin{myequation}
\label{eq:coend}
\coprod_{[n]\mapsto [m]} A^n \tens K_m \dbra \coprod_{[n]} A^n
\tens K_n
\end{myequation}
where $A^n\tens K_m$ is shorthand for a coproduct of copies of $A^n$
indexed by the set $K_m$.  There are adjunctions
\begin{myequation}
\label{eq:adjunct}
 \cM(A\tens K,X) \iso c\cM(A,X^K) \iso \sSet(K,\cM(A,X)) 
\end{myequation}
for $A\in c\cM$, $K\in \sSet$, and $X\in \cM$.  Here $X^K$ is the
cosimplicial object $[n]\mapsto X^{K_n}$, where $X^{K_n}$ denotes a
product of copies of $X$ indexed by the set $K_n$.
 One checks---using the above adjunctions or otherwise---that $A\tens
\del{n}\iso A^n$, and $A\tens \bd\del{n}$ is isomorphic to the $n$th
latching object of $A$ \cite[Def. 15.2.5]{H}.  See \cite[Section
4]{drep} for the dual situation with $s\cM$ instead of $c\cM$.

Write $s^2\Set$ for the category of bisimplicial sets and $c^2\cM$ for
the category of bi-cosimplicial objects in $\cM$.  When drawing a
bisimplicial set $P$ we will draw each $P_{m,*}$ horizontally, and
each $P_{*,n}$ vertically.  If $K\in \sSet$ and $P\in s^2\Set$,
let $\vMap(K,P)$ denote the simplicial set $[n]\mapsto
\sSet(K,P_{*,n})$.  We are mapping $K$ into the {\it vertical\/}
simplicial sets of $P$.  
 
If $K,L\in \sSet$ write $K\extens L$ for
the bisimplicial set $[m],[n] \mapsto K_m\times L_n$.  Observe that
there is an adjunction formula
\begin{myequation}
\label{eq:adjunct2}
 s^2\Set(K\extens L,P)\iso \sSet(L,\vMap(K,P)).
\end{myequation}
Note in particular that $s^2\Set(\del{m}\extens \del{n},P)\iso
P_{m,n}$.  

If $P\in s^2\Set$ and $A\in c^2\cM$, one can form a coend $A\tens P
\in \cM$ similarly to what was done in (\ref{eq:coend}).  There are
adjunction formulas analogous to (\ref{eq:adjunct}).  One checks that
$A\tens (\del{m}\extens \del{n})\iso A^{m,n}$, and more generally
$A\tens (\del{m}\extens L)\iso A^{m,*}\tens L$ (use
(\ref{eq:adjunct2}) for both).  So, for instance,
$A\tens (\del{m}\extens \bd\del{n})$ is the $n$th latching object for
the cosimplicial object $A^{m,*}$.

Finally, recall from \cite[p. 125]{BF} that the diagonal functor $\diag\colon
s^2\Set \ra \sSet$ has a left adjoint which we will call $d\colon
\sSet \ra s^2\Set$.  It follows immediately from adjointness that
$d\del{n} \iso \del{n}\extens\del{n}$.  Since $d$ preserves colimits
and every simplicial set is a colimit of $\del{n}$'s, this tells us
what $d$ does to any simpiclial set.

By chasing through adjunctions one finds that if
$X\in c^2\cM$ and $K\in \sSet$ then $\diag(X) \tens K \iso X\tens dK$.

\begin{proof}[Proof of Lemma~\ref{le:addlem}]
Consider the object $X\in c^2\cM$ given in the statement of the lemma.
Our task is to show that $\diag(X)\tens \bd\del{n} \ra \diag(X)\tens
\del{n}$ is a cofibration, for each $n$.  This is the condition for
$\diag(X)$ to be Reedy cofibrant.  Using the isomorphisms
$\diag(X)\tens K \iso X\tens dK$, this is equivalent to showing that
the map $X\tens d(\bdd{n}) \ra X\tens d\del{n}\iso X\tens (\del{n}\extens
\del{n})$ is a cofibration.

Let $S$ denote the set of all maps $P\ra Q$ of bisimplicial sets such
that $X\tens P \ra X\tens Q$ is a cofibration in $\cM$.  This set is
closed under composition and cobase change.  Our assumption about the
latching maps of $X$ amounts to saying that the maps
\begin{myequation}
\label{eq:boxprod}
 (\bdd{k}\extens \del{n}) \amalg_{(\bdd{k}\extens \bdd{n})}
(\del{k}\extens \bdd{n}) \lra \del{k}\extens \del{n}
\end{myequation}
belong to $S$, for all $n$ and $k$.  Our goal is to show that this
forces $d(\bdd{n}) \ra \del{n}\extens \del{n}$ to also belong to $S$.

We have now reduced things to a problem in combinatorial homotopy
theory.  Namely, we must show that $d(\bdd{n}) \ra \del{n}\extens
\del{n}$ can be obtained from the maps in (\ref{eq:boxprod}) by
iterated cobase changes and compositions.  But a little thought shows
that {\it every\/} monomorphism of bisimplicial sets can be obtained
in this way from the maps in (\ref{eq:boxprod}) (the point is that
every monomorphism of simplicial sets can be obtained from the maps
$\bdd{n} \ra\del{n}$ in the same way).  So we are done.
\end{proof}


\section{Universal additive model categories}
\label{se:univ}

Suppose $\cC$ is a small category.  The paper \cite{duniv} introduced
the idea of a universal model category built from $\cC$, there denoted
$U\cC$.  This is just the category of functors $\Func(\cC^{op},\sSet)$
with a well-known model structure.  

If $\cC$ is also an additive category then one can ask for a
universal {\it additive\/} model category built from $\cC$.  This
section develops something along these lines, although the `universal'
properties are slightly weaker than one might hope for.  They are
enough for reproducing the enrichment results of \cite{hend}, however.

\medskip  

\subsection{Presheaves and additive presheaves}
Let $\cC$ be a small, additive category.  Let $\Func(\cC^{op},\Ab)$
denote the category of all functors.  Note that for every $X\in \cC$,
the representable functor $rX\colon \cC^{op}\ra \Ab$ defined by
$U\mapsto \cC(U,X)$ is additive.

The Yoneda Lemma does not hold in $\Func(\cC^{op},\Ab)$: that is, if
$F\in \Func(\cC^{op},\Ab)$ one need not have $\Hom(rX,F)\iso F(X)$ for
all $X\in \cC$.  But it is easy to check that this {\it does\/} hold
when $F$ is an additive functor.

Let $\Func_{ad}(\cC^{op},\Ab)$ denote the full subcategory of additive
functors.  The following lemma records several basic facts about this
category.

\begin{lemma} 
\label{le:functorcats}
Let $\cC$ be a small, additive category.
\begin{enumerate}[(a)]
\item Colimits and limits in $\Func_{ad}(\cC^{op},\Ab)$ are the same
as those in $\Func(\cC^{op},\Ab)$.  
\item Every additive functor $F\in \Func(\cC^{op},\Ab)$ is isomorphic
to its canonical colimit with respect to the embedding $r\colon \cC
\inc \Func(\cC^{op},\Ab)$.  That is, the natural map
$\Bigl [ \colim\limits_{rX \ra F} (rX) \Bigr ] \ra F$ is an isomorphism.
\item
The additive functors in $\Func(\cC^{op},\Ab)$ are precisely those
functors which are colimits of representables.
\item The inclusion $i\colon \Func_{ad}(\cC^{op},\Ab) \inc
\Func(\cC^{op},\Ab)$ has a left adjoint $\Ad$ (for `additivization'),
and the composite $\Ad \circ i$ is naturally isomorphic to the
identity.
\item Suppose given a co-complete, additive category $\cA$ and an
additive functor $\gamma\colon \cC \ra \cA$.  Define $\Sing\colon \cA \ra
\Func_{ad}(\cC^{op},\Ab)$ by letting $\Sing(a)$ be the
functor $c\mapsto \cA(\gamma c,a)$.  Then $\Sing$ has a left adjoint
$\Real$, and there are natural isomorphisms $\Real(rX) \iso
\gamma(X)$.
\end{enumerate}
\end{lemma}

\begin{proof}
We mostly leave this to the reader.  We note, however, that the
fact that $\cC$ has finite coproducts (which is part of the
definition of an additive category) is needed in (b).  This
ensures that the categories indexing the canonical colimits are
pseudo-filtered, in the sense that for any objects $i$ and $j$ there
is a third object $k$ and maps $i\ra k$, $j\ra k$.  

Also, we define the additivization functor from (d).  If $F$ is any
functor, then $(\Ad F)(X)$ is the quotient of $F(X)$ by the subgroup
generated by all $(f+g)^*(s)-f^*(s)-g^*(s)$ for all objects $Y$, all
functions $f,g\colon X\ra Y$, and all $s\in F(Y)$.  

For the proof of (e), note the following.  If $X\in \cA$ and $B$ is an
abelian group, one can define $X\tens B$ as the coequalizer of two
maps
\[ \coprod_{B\times B} X \dbra \coprod_B X.\]
Here the objects are coproducts of copies of $X$, indexed by the sets
$B\times B$ and $B$, respectively.  To describe the two maps, we have
to say what they do to each summand corresponding to a pair
$(b_1,b_2)$.  The first map is just the inclusion into the summand
indexed by $b_1+b_2$.  The second map is the sum of the two inclusion maps
corresponding to the summands $b_1$ and $b_2$.  
One checks that with this definition there is a natural adjunction
isomorphism $\cA(X\tens B,Y)\iso \Ab(B,\cA(X,Y))$.

Recall that we are given an additive functor $\gamma\colon \cC\ra \cA$.
Given a functor $F\colon \cC^{op}\ra \Ab$, we consider the coend
\[ \gamma\tens F=\coeq \Biggl [ \coprod_{c\ra d} \gamma(c) \tens F(d)
\dbra \coprod_c \gamma(c)\tens F(c) \Biggr ].
\]
When $F$ is an additive functor one defines $\Real(F)=\gamma\tens F$.  
It is routine to check that this is a left adjoint to $\Sing$.
\end{proof}

By \cite[Th. 11.6.1]{H} the category $\Func(\cC^{op},s\Ab)$ has a
cofibrantly-generated model structure in which the weak equivalences
and fibrations are defined objectwise.  We will need the analogous
result for the category of additive functors:

\begin{lemma}
\label{le:functor-mc}
Let $\cC$ be a small, preadditive category.
Then $\Func_{ad}(\cC^{op},s\Ab)$ has a cofibrantly-generated
model structure in which the weak equivalences and fibrations are
defined objectwise.  This model structure is simplicial, left proper,
combinatorial, and cellular.
\end{lemma}

\begin{proof}
The proof uses the adjoint pair $(\Ad,i)$ to create the model
structure, as in \cite[Th. 11.3.2]{H}.  Recall that the model category
$\Func(\cC^{op},s\Ab)$ has generating trivial cofibrations $J=\{rX
\times \Z[\Lambda^{n,k}] \ra rX \times \Z[\Delta^n] \mid X\in \cC\}$.
Our
notation is that if $K\in \sSet$ then $\Z[K] \in s\Ab$ is the
levelwise free abelian group on $K$; and if $A\in s\Ab$ then $rX\times
A$ denotes the presheaf $U \mapsto \cC(U,X) \times A$ (with the product
performed levelwise).  Note that we think of $rX$ as a $\Set$-valued
functor here, so $\cC(U,X)\times A$ denotes a direct sum of copies of
$A$ indexed by the set $\cC(U,X)$---this is {\it not\/} the same as as
the direct product of the abelian groups $\cC(U,X)$ and $A$.

To apply \cite[11.3.2]{H} we must verify that the functor $i$ takes
relative $\Ad(J)$-cell complexes to weak equivalences.  However, note
that if $A$ is an abelian group then $\Ad(rX\times A)\iso rX\tens A$,
where the latter refers to the presheaf $U\mapsto \cC(U,X)\tens A$.
So $\Ad(J)$ is the set of maps $rX\tens \Z[\Lambda^{n,k}] \ra rX
\tens \Z[\Delta^n]$.  Objectwise, these maps are  monomorphisms and
 weak equivalences of simplicial abelian groups.

Now, the model category $s\Ab$ has the special property that a pushout
of a map which is both a monomorphism and a weak equivalence is still
a monomorphism and weak equivalence.  The fact that forming pushouts
in $\Func_{ad}(\cC^{op},s\Ab)$ and $\Func(\cC^{op},s\Ab)$ give the
same answers (by Lemma~\ref{le:functorcats}(a)) and are done
objectwise therefore shows that the $\Ad(J)$-cell complexes are
objectwise monomorphisms and objectwise weak equivalences.  In
particular, they are weak equivalences in $\Func(\cC^{op},s\Ab)$.

Finally, it is routine to check that the resulting model structure is
simplicial, left proper, combinatorial, and cellular.
\end{proof}

From now on we will write $\Ua \cC$ for the category
$\Func_{ad}(\cC^{op},s\Ab)$ with the model structure provided by the
above lemma.  The reason for the notation is provided by the next result.

Recall that if $L_1,L_2\colon \cM \ra \cN$ are two Quillen maps then a
\dfn{Quillen homotopy} from $L_1$ to $L_2$ is a natural transformation
$L_1 \ra L_2$ which is a weak equivalence on the cofibrant objects.

If $\cM$ is a  model category and $S$ is a set of maps in $\cM$,
then we use $\cM/S$ to denote the left Bousfield localization of $\cM$
at $S$, if it exists.   See \cite[Chapters 3--4]{H} and \cite{duniv}
for a discussion.  The localizations always exist when $\cM=\Ua\cC$,
since this model category is left proper and cellular.

\begin{thm} 
\label{th:additive-univ}
Let $\cM$ be an additive model category.
\begin{enumerate}[(a)]
\item Let $\cC$ be a small, additive category and $\gamma\colon
\cC \ra \cM$ an additive functor taking values in the cofibrant
objects.  Then there is a Quillen pair $\Real \colon \Ua\cC \adjoint
\cM \colon \Sing$ together with a natural weak equivalence $\Real
\circ r \we \gamma$.  Moreover, any two such Quillen pairs are
connected by a zig-zag of Quillen homotopies.
\item If $\cM$ is combinatorial then there
is a Quillen equivalence $\Ua\cC/S \we \cM$ for some small, additive
category $\cC$ and some set of maps $S$ in $\Ua\cC$.
\item Suppose $\cM \bwe \cM_1 \we \cdots \bwe \cM_n \we \cN$ is a
zig-zag of Quillen equivalences in which all the model categories are
additive.  If $\cM$ is combinatorial, there is a simple zig-zag of
equivalences 
\[ \cM \bwe \Ua\cC/S \we \cN 
\]
such that the derived equivalence $\ho(\cM)\he \ho(\cN)$ is
isomorphic to the derived equivalence given by the original zig-zag.
\end{enumerate}
\end{thm}

\begin{proof}
For (a), one shows that giving a Quillen pair $\Real\colon \Ua\cC
\adjoint \cM \colon \Sing$ together with a natural weak equivalence
$\Real(rX)\we \gamma(X)$  is precisely the same as giving an additive
cosimplicial resolution on $\gamma$.  The proof of this is exactly the
same as \cite[Prop. 3.4]{duniv}.  Giving a Quillen homotopy between
two such Quillen pairs exactly amounts to giving a natural weak
equivalence between the corresponding cosimplicial resolutions.  This
proves (a), once one  recalls our definition of additive model categories.

The proof for (c) now exactly follow the case for $\U\cC$ given in
\cite[Cor. 6.5]{duniv}.  One uses along the way that adjoint functors
between additive categories are necessarily additive functors.

The proof of (b) is slightly more complicated; we will return to it
at the end of this section, after some discussion.
\end{proof}

\begin{remark}
\label{re:counter}
The result in (c) is false if one does not assume that all the
$\cM_i$'s are additive.  For an example, let $R$ be the dga $\Z[e;
de=2]/(e^4)$ and let $T$ be the dga $\Z/2[x;dx=0]/(x^2)$, where $e$
has degree $1$ and $x$ has degree $2$.  Let $\cM$ and $\cN$ be the
categories of $R$- and $T$-modules, respectively.  These turn out to
be Quillen equivalent, but they cannot be linked by a zig-zag of
Quillen equivalences between {\it additive\/} model categories.  A
verification of these claims can be found in \cite[Section 8]{tpwe}.
\end{remark}

\subsection{Additive presentations}
We turn to the proof of Theorem~\ref{th:additive-univ}(b).  This will
be deduced from the work of \cite{comb} plus some purely formal
considerations.

Let $\cM$ be a combinatorial model category.  By \cite[Prop. 3.3]{comb},
there is a small category $\cC$ and a functor $\cC \ra \cM$ such that
the induced map $L\colon \U\cC \ra \cM$ is {\it homotopically
surjective\/} (see \cite[Def. 3.1]{comb} for the definition).  Then
\cite[Prop. 3.2]{comb} shows that this fact implies there is a set of
maps $S$ in $\U\cC$ which the derived functor of $L$ takes to weak
equivalences, and such that the resulting map $\U\cC/S \ra \cM$ is a
Quillen equivalence.

Now suppose that $\cM$ was also an additive model category.  By
examining the proof of \cite[Prop. 3.3]{comb} one sees that the $\cC$
constructed there is actually an additive category and the functor
$\gamma\colon\cC \ra \cM$ an additive functor taking values in the
cofibrant objects (the category $\cC$ is a
certain full subcategory of the cosimplicial objects over $\cM$).  By
Theorem~\ref{th:additive-univ}(a) there is an induced map
$F\colon\Ua\cC \ra \cM$.  Again using \cite[Prop. 3.2]{comb}, it will be
enough to prove that this map is homotopically surjective.

Consider now the following sequence of adjoint pairs:
\[ \xymatrix{
\Func(\cC^{op},\sSet) \ar@<0.5ex>[r]^{\Z}
& \Func(\cC^{op},s\Ab) \ar@<0.5ex>[r]^{\Ad}\ar@<0.5ex>[l]^U
& \Func_{ad}(\cC^{op},s\Ab) \ar@<0.5ex>[r]^-F\ar@<0.5ex>[l]^i
& \cM \ar@<0.5ex>[l]^-\Sing
}
\]
The composite of the right adjoints is clearly the right adjoint of
$L$, so the composite of the left adjoints is $L$.  We have
constructed things so that this composite is homotopically surjective,
and we are trying to show that $F$ is also homotopically surjective.

In the following lemma, note that the presheaf $rX$ can be regarded as
an object of either $\Func_{ad}(\cC^{op},\Ab)$ or
$\Func(\cC^{op},\Set)$.  It will usually be clear from context which
one we intend.

\begin{lemma}
\label{le:AdZ}
If $X\in \cC$ then $\Ad(\Z (rX))\iso rX$. 
Said equivalently, one has $\Ad(\Z (Ui(rX))) \iso rX$.
\end{lemma}

\begin{proof}
This is clear, since the two functors $\Func_{ad}(\cC^{op},\Ab) \ra \Ab$
given by $F\mapsto \Func_{ad}(\Ad(\Z(rX)),F)$ and $F\mapsto
\Func_{ad}(rX,F)$ are both naturally isomorphic to $F\mapsto F(X)$.
\end{proof}

Let $G \in \Func_{ad}(\cC^{op},s\Ab)$.  Let $QG$ be the simplicial
presheaf whose $n$th level is
\[ \coprod_{rX_n \ra rX_{n-1} \ra \cdots \ra rX_0 \ra G_n} (rX_n)
\]
where the coproduct is in $\Func(\cC^{op},\sSet)$.  The simplicial
presheaf $QG$ is treated in detail in \cite[Sec. 2.6]{duniv}, as it is a
cofibrant-replacement functor for $\U\cC$.   Likewise,
let $Q_{ad}G$ be the simplicial presheaf whose $n$th level is
\[ \bigoplus_{rX_n \ra rX_{n-1} \ra \cdots \ra rX_0 \ra G_n} (rX_n)
\]
where the coproduct is now in $\Func(\cC^{op},s\Ab)$.  The proof of
\cite[Prop. 2.8]{comb} showing that $Q$ is a cofibrant-replacement
functor for $\U\cC$ adapts verbatim to show that $Q_{ad}$ is a
cofibrant-replacement functor for $\Ua\cC$.  Note that by
Lemma~\ref{le:AdZ} we have $Q_{ad}G\iso\Ad(\Z (Q(UiG)))$, since
$\Ad$ and $\Z(\blank)$ are left adjoints and therefore
preserve coproducts.  

Finally we are in a position to conclude the

\begin{proof}[Proof of Theorem~\ref{th:additive-univ}(b)]
We have reduced to showing that $F\colon \Ua\cC \ra\cM$ is
homotopically surjective.  Let $\Sing$ be the right adjoint of $F$.
Then we must show that for every fibrant object $X\in \cM$ the induced
map $FQ_{ad}(\Sing X)\ra X$ is a weak equivalence.

However, we have seen above that
\[ F[Q_{ad}\Sing X] \iso F[\Ad \Z (Q Ui(\Sing X))] \iso L[QUi(\Sing
X)]. 
\]
Recall that $Ui\Sing$ is the right adjoint to $L$.  
Since $L\colon \U\cC \ra \cM$ is homotopically
surjective we know  $LQ (Ui\Sing X) \ra X$ is a weak equivalence in
$\cM$, so we are done.
\end{proof}


\section{Homotopy enrichments over $\Spe^\Sigma(s\Ab)$}
\label{se:additive}

In this section and the next we prove the main results stated in
Section~\ref{se:intro}.  Except for the work in the next section, the
proofs are essentially the same as in \cite{hend}---but they use
Theorem~\ref{th:additive-univ} in place of \cite[Prop. 5.5]{hend}.  

\medskip

\subsection{Background on ring objects}
\label{se:background}
If $\cM$ is a monoidal model category which is combinatorial and
satisfies the monoid axiom, then by \cite[Th. 4.1(3)]{SS1} the
category of monoids in $\cM$ has an induced model structure where the
weak equivalences and fibrations are the same as those in $\cM$.
We'll write $\Ring[\cM]$ for this model category.  If $\cN$ is another
such monoidal model category and $L\colon \cM \adjoint \cN \colon R$
is a Quillen pair which is weak monoidal in the sense of
\cite[Def. 3.6]{SS3}, then there is an induced Quillen map
$\Ring[\cM] \ra \Ring[\cN]$.  This is a Quillen equivalence if $\cM
\ra \cN$ was a Quillen equivalence and the units in $\cM$ and $\cN$
are cofibrant \cite[Th. 3.12]{SS3}.

The adjunction $\Set_* \adjoint \Ab$ is strong monoidal, and therefore
induces strong monoidal Quillen functors $\Spe^\Sigma(\sSet_*)\adjoint
\Spe^\Sigma(s\Ab)$.  Therefore one gets a Quillen pair $F\colon
\Ring[\Spe^\Sigma] \adjoint \Ring[\Spe^\Sigma(s\Ab)]\colon U$.  By the
\dfn{Eilenberg-Mac\,Lane ring spectrum} associated to an $R\in
\Ring[\Spe^\Sigma(s\Ab)]$ we simply mean the ring spectrum $UR$.

\subsection{Additive enrichments}
Let $\cM$ be an additive, stable, combinatorial model category.  By
Theorem~\ref{th:additive-univ} there is a Quillen equivalence $\Ua
\cC/S \ra \cM$ for some small, additive category $\cC$ and some set of
maps $S$ in $\Ua\cC$.  The category $\Ua\cC/S$ is simplicial, left
proper, and cellular, so using \cite[Sections 8, 9]{stab} we may form
$\Spe^\Sigma (\Ua\cC/S)$.  Since $\Ua\cC/S$ is stable (since $\cM$
was), we obtain a zig-zag of Quillen equivalences
\[ \cM \bwe \Ua\cC/S  \we \Spe^\Sigma(\Ua\cC/S).
\]
Applying $\ME_0(\blank,\Sp^\Sigma(s\Ab))$ to this zig-zag gives a
diagram of bijections by \cite[3.14(d)]{hend}.

The category $\Ua\cC$ is a $s\Ab$-model category, and therefore
$\Spe^\Sigma(\Ua\cC/S)$ is a $\Spe^\Sigma(s\Ab)$-model category by
\cite[8.3]{stab}.  So $\Spe^\Sigma(\Ua\cC/S)$ comes with a natural
model enrichment by $\Spe^\Sigma(s\Ab)$, as in \cite[Ex. 3.2]{hend}.  We can
transport this enrichment onto $\cM$ via the Quillen equivalences, and
therefore get an element $\sigma_\cM\in \ME_0(\cM,\Spe^\Sigma(s\Ab))$.
Just as in \cite[Prop. 6.1]{hend}, one shows (using
Theorem~\ref{th:additive-univ}) that this quasi-equivalence class does
not depend on the choice of $\cC$, $S$, or the Quillen equivalence
$\Ua\cC/S \we \cM$.

We can now give the:

\begin{proof}[Proof of Theorem~\ref{th:additive-enrich}]
We have just constructed the enrichment $\sigma_\cM$.  The proof that it
is preserved by Quillen equivalences is exactly the same as in
\cite[Prop. 6.2]{hend}, but using Theorem~\ref{th:additive-univ}.
\end{proof}

Let $X\in \cM$, and let $\tilde{X}$ be a cofibrant-fibrant object
weakly equivalent to $X$.  We write \mdfn{$\hdga(X)$} for any object
in $\Ring[\Spe^\Sigma(s\Ab)]$ having the homotopy type of
$\sigma_\cM(\tilde{X},\tilde{X})$, and we'll call this the \mdfn{additive
homotopy endomorphism object of $X$}.  By
\cite[Cors. 3.6, 3.7]{hend} this
homotopy type depends only on the homotopy type of $X$ and the
quasi-equivalence class of $\sigma_\cM$---and so it is a well-defined
invariant of $X$ and $\cM$.

\begin{proof}[Proof of Proposition~\ref{pr:hdga-invariance}]
This is entirely similar to the proof of \cite[Th. 1.4]{hend}, but
using Theorem~\ref{th:additive-univ}(c).
\end{proof}

\begin{proof}[Proof of Proposition~\ref{pr:hdga-spectral}]
Same as the proof of \cite[Prop. 1.5]{hend}.
\end{proof}

\begin{proof}[Proof of Proposition~\ref{pr:hdga=EM}]
We know that there exists a zig-zag of Quillen equivalences $\cM \bwe
\Ua\cC/S \we \Spe^\Sigma(\Ua\cC/S)$.  Therefore, using
\cite[Thm. 1.4]{hend} and Proposition~\ref{pr:hdga-invariance}
we may as well assume $\cM=\Spe^\Sigma(\Ua \cC/S)$.  This is an
$\Spe^\Sigma(s\Ab)$-model category, and so for any object $X$ we have a
ring object $\und{\cM}(X,X)$ in $\Spe^\Sigma(s\Ab)$.  The adjoint functors
$\Set_* \adjoint \Ab$ induce a strong monoidal adjunction $F\colon
\Spe^\Sigma(\sSet_*) \adjoint \Spe^\Sigma(s\Ab)\colon U$.  The
$\Spe^\Sigma(s\Ab)$-structure on $\cM$ therefore yields an induced
$\Spe^\Sigma$-structure as well (see \cite[Lem. A.5]{hend}).  In
this structure, the endomorphism ring spectrum of $X$ is precisely
$U[\und{\cM}(X,X)]$.  Using \cite[Prop. 1.5]{hend}, we know that this has
the homotopy type of the ring spectrum $\hEnd(X)$, at least when $X$
is cofibrant-fibrant.  And Proposition~\ref{pr:hdga-spectral} says
that $\und{\cM}(X,X)$ has the homotopy type of $\hdga(X)$.  This is all we
needed to check.
\end{proof}


\section{Chain enrichments}
\label{se:chain}
 
Proposition~\ref{pr:hdga-spectral} says that if $\cM$ is a
$\Sp^\Sigma(s\Ab)$-model category then one can compute $\hdga(X)$
using the $\Sp^\Sigma(s\Ab)$-structure.  We would like to prove a
similar result for $\ChZ$-model categories, where $\ChZ$ denotes the
model category of unbounded chain complexes of abelian groups.  These are
what arise most commonly in algebraic situations.

The monoidal model categories $\Sp^{\Sigma}(s\Ab)$ and $\ChZ$ can be
connected by a zig-zag of weak monoidal Quillen equivalences, as
described in \cite{S}.  This zig-zag can be used to translate
enrichment-type information between these two categories.  However,
this is not as straightforward as one might expect; there are
complications arising from the monoidal properties of the Dold-Kan
equivalence between $s\Ab$ and $\ch_+$, as analyzed in~\cite{SS3}.
Our method for dealing with this requires some cumbersome machinery
and gives a slightly weaker result than one would like.  However, it
is the best we can do at the moment.

\subsection{Statement of the result}
We give $\ChZ$ the projective model structure, where weak equivalences
are quasi-isomorphisms and fibrations are surjections.  Recall again
that a \mdfn{$\ChZ$-model category} is a model category with
compatible tensors, cotensors, and enrichments over $\ChZ$
satisfying an analogue of SM7; see Section~\ref{se:enrich}.
For $X, Y$ in $\cM$,
we denote the enriched hom-object in $\ChZ$ by $\Fhomu(X,Y)$.

Note that a $\ChZ$-model category is automatically additive
and stable.  See Corollary~\ref{co:Cmodel=additive} for the
additivity, 
and~\cite[3.5.2]{SS2} 
or~\cite[3.2]{GS} for stability. 

Recall from \cite{S} that
there are two Quillen equivalences
\[
\xymatrix{
\Sp^\Sigma(\ch_+) \ar@<0.5ex>[r]^-D & \Ch \ar@<0.5ex>[l]^-R 
& \text{and} &
   \Sp^\Sigma(\ch_+) \ar@<0.5ex>[r]^-L &\Sp^\Sigma(s\Ab) 
\ar@<0.5ex>[l]^-{\namma}
}
\]
in which $(D,R)$ is strong monoidal and $(L, \namma)$ is weak monoidal.
These induce Quillen equivalences between the corresponding
model categories of rings:
\begin{myequation}
\label{eq:rings}
 \Ring(\Sp^\Sigma(\ch_+)) \ladjoint{\sim} \DGA \qquad
\Ring(\Sp^\Sigma(\ch_+)) \ladjoint{\sim} \Ring(\Sp^\Sigma(s\Ab)).
\end{myequation}
In the first equivalence of (\ref{eq:rings}) the left and right adjoints 
are just the restrictions of $D$ and $R$, as these were strong monoidal.  
In the second, the right adjoint is just $\namma$ again,
but the left adjoint is more complicated; see~\cite[3.3]{SS3}.

Let $\und{\namma}$ and
$\und{D}$ denote the derived functors of $\namma$ and $D$ from
(\ref{eq:rings}), and write $\Theta'=\und{D}\,\und{\namma}$.  So
$\Theta'$ is a functor
\[ \Ho(\Ring(\Sp^\Sigma(s\Ab))) \ra \Ho(\DGA).\]

Let $\cM$ be a stable, combinatorial, additive model category and let
$X\in \cM$.  We have shown how to associate to $X$ an object
$\hEnd_{ad}(X) \in \Ring(\Sp^\Sigma(sAb))$.  By applying $\Theta'$ we
get the \dfn{homotopy endomorphism dga} of $X$.  Denote this as
$\hdgadga(X)=\Theta' \bigl [ \hEnd_{ad}(X) \bigr ]$.

The goal for this section is to prove 
Proposition~\ref{pr:hdga-chain}.  We restate the result here for the
convenience of the reader.

\begin{prop}
\label{pr:hendCH}
Suppose that $\cM$ is a combinatorial $\ChZ$-model category,
and that $\cM$ has a generating set of compact objects.  Let $X\in
\cM$ be cofibrant and fibrant.  Then the dga $\Fhomu(X,X)$ is
quasi-isomorphic to $\hdgadga(X)$ 
\end{prop}

Proposition~\ref{pr:hendCH} will be proven by reducing from a
$\ChZ$-model category to a $\Sp^\Sigma(s\Ab)$-model category and then
applying results of Section~\ref{se:additive}. The reduction from
$\ChZ$ to $\Spes(\ch_+)$ will be simple because of the strong monoidal
equivalence between these two categories.  The following proposition
provides the reduction from $\Spes(\ch_+)$ to $\Spes(s\Ab)$.  This is
where all the enriched category theory from 
Sections~\ref{se:enrich} through \ref{se:transfer} is needed.  Recall
that for a general $\cD$-model category $\cN$ we denote the morphism
object in $\cD$ by $\ND(X,Y)$.

\begin{prop}
\label{pr:ch->sAb}
Let $\cM$ be a combinatorial  $\Spes(\ch_+)$-model category with a
generating set of compact objects.  Let $X\in \cM$ be a
cofibrant-fibrant object.  Then there exists
\begin{enumerate}[(i)]
\item a combinatorial, $\Sp^\Sigma(s\Ab)$-model category $\cN$,
\item a zig-zag of Quillen equivalences between $\cM$ and $\cN$, where
the intermediate model categories are all additive, and
\item a cofibrant-fibrant object $Y\in \cN$
\end{enumerate}
such that $Y$ is taken to $X$ by the derived functors of the Quillen
equivalences and $\und{\nu} \bigl [ \Nsab(Y,Y) \bigr ]$ is weakly equivalent 
to $\und{\cM}_{\Spes(\ch_+)}(X,X)$. 
\end{prop}

\begin{proof}
This is a special case of Proposition~\ref{pr:technical} and
Corollary~\ref{cor:technical}.
We need to verify the properties for
$\cC=\Sp^{\Sigma}(\ch_+)$ and $\cD=\Sp^{\Sigma}(s\Ab)$ stated just
prior to Proposition~\ref{pr:technical},
with $(F,G)$
replaced by $(L,\nu)$.  Axioms (QI1-2) for $\cC$ follow from
\cite[3.2, 3.3]{S}.  The fact that $(L,\nu)$ is a weak monoidal Quillen
equivalence is given in \cite[4.3]{S}.  All the other conditions are
easy exercises, but see also Example~\ref{ex1} for more information.
\end{proof}

Using the above proposition, we can complete the following:

\begin{proof}[Proof of Proposition~\ref{pr:hendCH}]
Let $\cM$ be a combinatorial  $\Ch$-model category with a
generating set of compact objects.  Let $\cC= \Spes(\ch_+)$.  Using
the strong monoidal adjunction $(D,R)$, $\cM$ becomes a
$\cC$-model category via the definitions $Z\tens
c=Z\tens D(c)$, $Z^c=Z^{Dc}$, and $\mM_{\cC}(W,Z)=R[\mM_{\Ch}(W,Z)]$
where $W,Z\in \cM$ and $c\in \cC$.  See \cite[Lem. A.5]{hend}.

Now we apply Proposition~\ref{pr:ch->sAb} to $\cM$ with this
$\cC$-model structure to construct $\cN$ and $Y$.  By
Proposition~\ref{pr:hdga-invariance}, the additive homotopy
endomorphism spectra corresponding to $X$ and $Y$ are weakly
equivalent.  Let $\cD =\Spes(s\Ab)$.  Since $\cN$ is a
$\cD$-model category, we have by
Proposition~\ref{pr:hdga-spectral} that $\hEnd_{ad}(Y)$ is weakly
equivalent to $\und{\cN}_{\cD}(Y,Y)$.  So we have
\[
\hdgadga(X)=
\Theta'[\hdga(X)] \he 
\Theta'[\hdga(Y)] \he 
\und{D}\,\und{\nu}\bigl [\mN_{\cD}(Y,Y)\bigr ]
\]
(recalling that $\Theta'=\und{D}\,\und{\nu}$).

But $\cN$ and $Y$ were chosen
in such a way that we have 
$\und{\nu}\bigl [ \und{\cN}_{\cD}(Y,Y) \bigr ] \he \und{\cM}_{\cC}(X,X)$.
So in fact
\[
\hdgadga(X)\he \und{D}\,\bigl [ \und{\cM}_{\cC}(X,X)\bigr ]=\und{D} \,
\und{R}\bigl [ \mM_{\Ch}(X,X) \bigr ] \he \mM_{\Ch}(X,X).
\]
\end{proof}

\appendix
 \vspace{0.3in}

\section{Homotopy theory of $\cC I$-categories}\label{sec-ci-cat}
\label{se:CI-cat}

The present section reviews and expands on results from \cite{SS3}.  In
particular, \cite{SS3} often states results in settings which are
extremely general and therefore require somewhat awkward hypotheses.
Here we will specialize, replacing those hypotheses with conditions
more readily checked in practice.  

\medskip

We assume that $\cC$ is a combinatorial, symmetric monoidal model
category.  Also, $\cC$ is assumed to satisfy the {\em monoid axiom} of
\cite[3.3]{SS1}.  We'll refer to those conditions as our {\em
`standing assumptions'}.  Finally, we will sometimes require the
following two conditions as well:
\begin{itemize}
\item[(QI1)]
For any cofibrant object $A\in \cC$ and any weak
equivalence $X\ra Y$, the map $A\tens X \ra A\tens Y$ is also a weak
equivalence.
\item[(QI2)]
Suppose $A\cofib B$ is a cofibration, and $X$ is any object.  Then for
any map $A\tens X \ra Z$, the map from the homotopy pushout of $B\tens
X \lla A\tens X \ra Z$ to the pushout is a weak equivalence.
\end{itemize}
The abbreviation (QI) is for `Quillen invariance', as these conditions
will be used to check what \cite[3.11]{SS3} 
 calls Quillen invariance for modules.

\begin{example}
The category $\ch_+$ of non-negatively graded chain complexes with
tensor product and its usual `projective' model structure satisfies
(QI1-2).  It follows from \cite[3.2, 3.3]{S} 
 that
$\Sp^\Sigma(\ch_+)$ also satisfies (QI1-2).  Typically, these axioms
will follow from the existence of an `injective' model structure for
$\cM$ in which all objects are cofibrant, provided such a model
structure is a Quillen module over the corresponding projective
version.
\end{example}

Let $I$ be a set.  We assume the reader is familiar with the notion
of \mdfn{$\cC I$-category} (a category enriched over $\cC$ with object set
$I$) from \cite[6.2]{Bor}.  
If $\cO$ is a $\cC I$-category, then the category of right
$\cO$-modules (contravariant 
$\cC$-functors from $\cO$ to $\cC$) is defined in \cite[6.2]{Bor}; 
see also \cite[Section 6]{SS3}.

\begin{prop}
\label{pr:Omod}
Let $\cO$ be a $\cC I$-category.
\begin{enumerate}[(a)]
\item The category $\Mod\cO$ has a model category structure in which
the weak equivalences and fibrations are defined objectwise.
\item Let $\cO \ra \cR$ be a map of $\cC I$-categories.
Then there is a Quillen map $\Mod\cO \ra \Mod\cR$ in which the right
adjoint is restriction.  
\item
If $\cO \ra \cR$ is a weak equivalence and
$\cC$ satisfies (QI1-2), then the above Quillen map is a Quillen
equivalence.
\end{enumerate}
\end{prop}

\begin{proof}
Part (a) is \cite[6.1(1)]{SS3}. 
  For (b) we need only construct
the left adjoint, as the restriction clearly preserves fibrations and
trivial fibrations.  This construction is given in the paragraph above 
\cite[6.1]{SS3}. 
  Denote this left adjoint by $X\mapsto X\tens_{\cO} \cR$.

Part (c) requires a little work.  First, for any $i\in I$ and $A\in
\cC$ let $A\tens \Fr_i(\cO)$ be the `free $\cO$-module generated by
$A$ at spot $i$'.  This is defined by $j\mapsto A\tens \cO(j,i)$.  It
is easy to see that we have the adjunction $\Mod\cO(A\tens
\Fr_i(\cO),X)\iso \cC(A,X(i))$.  From this it immediately follows that
\[ [A\tens \Fr_i(\cO)] \tens_{\cO} \cR \iso A\tens \Fr_i(\cR).
\]
As a another consequence of the adjunction,
observe  that $\Mod\cO$ is cofibrantly-generated and the
generating cofibrations are maps of the form $A\tens \Fr_i(\cO) \ra
B\tens \Fr_i(\cO)$ where $A\cofib B$ is a generating cofibration of
$\cC$. 

By \cite[6.1(2)]{SS3}, 
 to prove (c) it suffices to check that for
any cofibrant $\cO$-module $N$ the natural map $N\ra U[N\tens_{\cO}
\cR]$ is a weak equivalence, where $U$ is the restriction $\Mod\cR \ra
\Mod\cO$.  Let $G$ denote the composite functor $X\mapsto U[X\tens_{\cO}
\cR]$, so that we are concerned with the natural transformation $\Id
\ra G$.  Note that when $X=A\tens \Fr_i(\cO)$ we have $G(X)=A\tens
\Fr_i(\cR)$.   If $A$ is cofibrant, the map $X\ra G(X)$ is an objectwise 
weak equivalence because $\cO \ra \cR$ is (this uses (QI1)).

Apply the small object argument to factor $\emptyset \ra N$ as
a cofibration followed by a trivial fibration.  This gives us a (possibly
transfinite) sequence of cofibrations
\[ \emptyset=W_0 \cofib W_1 \cofib W_2 \cofib \cdots \]
in which $W_{i+1}$ is obtained from $W_i$ by a pushout diagram
\[ \xymatrix{ \coprod_j A_j\tens \Fr_j(\cO) \ar[r]\ar[d] & W_i\ar[d] \\
  \coprod_j B_j \tens \Fr_j(\cO) \ar[r] & W_{i+1},
}
\]
together with a trivial fibration $W_\infty=\colim_i W_i \ra N$.
Since $N$ is cofibrant, $N$ is a retract of $W_\infty$.  So it
will suffice to show that $W_\infty\ra G(W_\infty)$ is a weak
equivalence, as $N\ra GN$ is a retract of this map.  

We first prove that if $W_{i-1} \ra G(W_{i-1})$ is a weak equivalence
then the same is true of $W_i\ra G(W_i)$.  To see this, note that
we have the following diagram:
\[ \xymatrix{ 
\coprod B_j\tens \Fr_j(\cO) \ar[d]^\sim & \coprod A_j \tens \Fr_j(\cO) 
\ar[d]^\sim \ar@{ >->}[l]\ar[r] &
W_{i-1} \ar[d]^\sim \\
\coprod G(B_j\tens \Fr_j(\cO)) & \coprod G(A_j\tens \Fr_j(\cO)) \ar[l]
\ar[r] &
G(W_{i-1}).
}
\]
The pushout of the top row is $W_i$, and of the bottom row is $G(W_i)$
(the latter follows because $G$ preserves colimits).  Note that
$G(A_j\tens \Fr_j(\cO)) \ra G(B_j\tens \Fr_j(\cO))$ is a cofibration,
as it is just the map $A_j\tens \Fr_j(\cR) \ra B_j\tens \Fr_j(\cR)$.
It follows that $G(W_{i-1})\ra G(W_i)$ is a cofibration.

Certainly the above diagram induces a weak equivalence of {\it homotopy\/}
pushouts.  We claim these homotopy pushouts are weakly equivalent to the
corresponding pushouts.  This is an objectwise question, since
pushouts, homotopy pushouts, and weak equivalences in the module
category are all determined objectwise.  The claim for the top row then
follows directly from (QI2).  The claim for the bottom row is
similar, but uses the identification $G(B_j\tens \Fr_j(\cO))=B_j\tens
\Fr_j(\cR)$, etc. 

Thus, we have shown that $W_i \ra G(W_i)$ is a weak equivalence
whenever $W_{i-1}\ra G(W_{i-1})$ is so.  It is trivial that $W_0\ra
G(W_0)$ is a weak equivalence.  
The result now follows by a transfinite
induction, using \cite[17.9.1]{H} to pass the weak equivalences to the
limit ordinals.  One again uses that $G$ preserves colimits.
\end{proof}

\begin{prop}  
\label{pr:CIcat}
Let $I$ be a fixed set, and let
$L\colon \cC \adjoint \cD\colon R$ be a weak monoidal Quillen pair
(see Section~\ref{se:monoidalfunc})
where both $\cC$ and $\cD$ satisfy our standing assumptions.
\begin{enumerate}[(a)]
\item
The category $\cC I-\Cat$ (and likewise $\cD I-\Cat$) has a model
category structure in which weak equivalences and fibrations are
defined objectwise.  
\item
There is a Quillen map $\cC I-\Cat \ra \cD I-\Cat$ in
which the right adjoint is `apply $R$ objectwise'.  The left adjoint
will be denoted $L^{\cD I}$.
\item
Suppose $\unit_{\cC}$ and $\unit_{\cD}$ are cofibrant.  If $\cO$ is a
cofibrant $\cC I$-category then there are weak equivalences $L[
\cO(i,j)] \ra (L^{\cD I}\cO)(i,j)$ for every $i,j \in I$.  These are
adjoint to the maps provided by the adjunction unit $\cO \ra R (L^{\cD
I}\cO)$.
\item Suppose $(L,R)$ is a Quillen equivalence and $\unit_\cC$,
$\unit_\cD$ are cofibrant.
Then the induced Quillen map $\cC I-\Cat \ra \cD I-\Cat$ is also a
Quillen equivalence.
\end{enumerate}
\end{prop}

\begin{proof}
Part (a)  is \cite[6.3(1)]{SS3}. 
  For part (b) we argue as
follows.  Recall the category $\cC I-\Graph$ from \cite[6.1]{SS3}, 
and that this category comes equipped with a monoidal product
$\tens$.  A $\cC I$-category is precisely a monoid with respect to
this tensor product.  The existence of the desired left adjoint
follows from Lemma~\ref{le:adjoint-monoid} below.  As the right
adjoint obviously preserves fibrations and trivial fibrations, we have
a Quillen pair.

For part (c), note that there is a Quillen map $\cC I-\Graph \ra \cC
I-\Cat$ in which the right adjoint is the forgetful functor.  The
model structure on $\cC I-\Cat$ is `created' by these adjoint functors
from the cofibrantly-generated model structure on $\cC I-\Graph$.  
\cite[6.4(1)]{SS3} 
 proves the desired claim in the case $\cO$ is a
cell complex, but since any cofibrant object is a retract of a
cell complex one immediately obtains the more general statement. 

Finally, we prove (d).  Note that since the functor $\cD I-\Cat \ra \cC
I-\Cat$ is just `apply $R$ objectwise', a map of fibrant objects
$X \ra Y$ in $\cD I-\Cat$ is a weak equivalence if and only if $RX \ra
RY$ is a weak equivalence in $\cC I-\Cat$.  So by
Lemma~\ref{le:Qequiv} below, we only need to show that if $\cO$ is a
cofibrant $\cC I$-category and $L^{\cD I}\cO \we \cA$ is a fibrant
replacement in $\cD I-\Cat$, then $\cO \ra R\cA$ is a weak
equivalence.

Since weak equivalences are detected objectwise, we must check that
$\cO(i,j) \ra R[\cA(i,j)]$ is a weak equivalence for every $i,j\in I$.
But $\cO$ is cofibrant, so each $\cO(i,j)$ is cofibrant in $\cC$ (see
\cite[6.3(2)]{SS3} ---this uses that $\unit_{\cC}$ is 
cofibrant).  And since $\cA$ is fibrant, each $\cA(i,j)$ is fibrant.
Using the Quillen equivalence $(L,R)$, we are therefore reduced to
checking that $L[\cO(i,j)] \ra \cA(i,j)$ is a weak equivalence.  But
we are really looking at the composite
\[ L[\cO(i,j)] \ra [L^{\cD I}\cO](i,j) \ra \cA(i,j). \]
The second map was assumed to be a weak equivalence, and the first map
is a weak equivalence by part (c).  So we are done.
\end{proof}

In this proof we used the following two lemmas. 

\begin{lemma}
\label{le:adjoint-monoid}
Let $\cC$ be a monoidal category which is complete and co-complete.
Assume that for any $X\in \cC$ the functors $X\tens (\blank)$ and
$(\blank)\tens X$ preserve filtered colimits.  Then
\begin{enumerate}[(a)]
\item The category of monoids in $\cC$ is co-complete.
\item If $\cB$ is another monoidal category, and $L\colon
\cB \adjoint \cC\colon R$ is an adjunction where $R$ is weak monoidal,
then $R$ induces a functor $\cC-\Monoid \ra \cB-\Monoid$ and this
functor has a left adjoint.
\end{enumerate}
\end{lemma}

\begin{proof}
Let $T\colon \cC \ra \cC$ be the `free algebra' monad, where 
\[ T(X)=\unit
\amalg X \amalg (X\tens X) \amalg \cdots.
\]  
The monoids in $\cC$ are
precisely the $T$-algebras.
Our assumptions imply that
$T$ preserves filtered colimits, so \cite[4.3.6]{Bor} 
 implies that
$\cC-\Monoid$ is co-complete.

Part (b) is an immediate consequence of (a) and \cite[4.5.6]{Bor}. 
\end{proof}

\begin{lemma}\cite[1.3.16]{hovey}
\label{le:Qequiv}
Let $L\colon\cM \adjoint \cN\colon R$ be a Quillen pair.  Assume the
following two conditions hold:
\begin{enumerate}[(i)]
\item If $X$ and $Y$ are fibrant objects in $\cN$, a map $X\ra Y$ is a
weak equivalence if $RX \ra RY$ is a weak equivalence.
\item For every cofibrant object $A\in \cM$ and every fibrant
replacement $LA \we Z$ in $\cN$, the composite map $A\ra RLA \ra RZ$ is
a weak equivalence.
\end{enumerate}
Then $(L,R)$ is a Quillen equivalence.
\end{lemma}

Here is the final result we will need:

\begin{prop} 
\label{pr:CImod}
Again assume that $L\colon \cC \adjoint \cD\colon R$ is a weak
monoidal Quillen pair, where $\cC$ and $\cD$ satisfy our standing
assumptions.   Also assume that $\unit_{\cC}$ and $\unit_{\cD}$ are
cofibrant.  
\begin{enumerate}[(a)]
\item If $\cA$ is a $\cD I$-category then there is a Quillen map  
$\Mod R\cA \ra \Mod \cA$ in which the right adjoint is `apply $R$
objectwise'.

\item Let $\cO$ be a $\cC I$-category.  Then there is a Quillen map
$\Mod\cO \ra \Mod (L^{\cD I}\cO)$ in which the right adjoint is the
composition of `applying
$R$ objectwise, then restricting across $\cO \ra R(L^{\cD I}\cO)$'.  

\item If $L\colon \cC \adjoint \cD \colon R$ is a Quillen
equivalence and $\cO$ is a cofibrant $\cC I$-category, then $\Mod \cO
\ra \Mod (L^{\cD I}\cO)$ is also a Quillen equivalence.
\end{enumerate}
\end{prop}

\begin{proof}
If $X$ is in the functor category $\cD^I$, let $T_{\cA}X\in \cD^I$
be the functor $j\mapsto \coprod_j X(j)\tens \cA(\blank,j)$.  Note
that this is a monad in an obvious way, and that the
$T_{\cA}$-algebras are precisely the $\cA$-modules.  We have the
diagram of categories
\[ \xymatrix{\cC^I & \cD^I \ar[l]^R \\
   \Mod (R\cA) \ar[u] & \Mod\cA \ar[l]^R\ar[u]
}
\]
where the vertical maps are forgetful functors.
By \cite[4.5.6]{Bor} 
 the map $\Mod\cA \ra \Mod (R\cA)$ has a left
adjoint, since $\Mod \cA$ is cocomplete.  This clearly gives a Quillen
pair.

For (b) we use the composite of the two Quillen maps  
\[ \Mod \cO \ra \Mod (RL^{\cD I}\cO) \ra \Mod (L^{\cD I}\cO).
\]
The first is provided by Proposition~\ref{pr:Omod}(b), induced by the
map $\cO \ra RL^{\cD I}\cO$.  The second comes from (a) of the present
result. 

Finally, part (c) is just \cite[6.5(1)]{SS3}. 
\end{proof}

\bibliographystyle{amsalpha}

\begin{thebibliography}{JTTW}

\bibitem[AR]{AR} J. Adamek and J. Rosicky, \emph{Locally presentable
and accessible categories}, London Math. Society Lecture Note Series
{\bf 189}, Cambridge University Press, 1994.

\bibitem[Bo]{Bor} F. Borceux, {\em Handbook of categorical algebra 2:
Categories and structures}, Cambridge University Press, 1994.

\bibitem[BF]{BF} A. K. Bousfield and E. M. Friedlander, {\em Homotopy
theory of $\Gamma$-spaces, spectra, and bisimplicial sets}, Geometric
Applications of Homotopy Theory II, Lect. Notes. in Math., vol. 658,
Springer-Verlag, New York, 1978, pp. 80--130.


\bibitem[Dr]{Dr} V. Drinfeld, {\em Quotients of dg-categories\/},
J. Algebra {\bf 272} (2004), no. 2, 643--691.

\bibitem[D1]{drep} D. Dugger, \emph{Replacing model categories with
simplicial ones}, Trans. Amer. Math. Soc. {\bf 353}, no. 12 (2001),
5003--5027.


\bibitem[D2]{duniv} D. Dugger, \emph{Universal homotopy theories}, 
Adv. Math. {\bf 164}, no. 1 (2001), 144--176.

\bibitem[D3]{comb} D. Dugger, {\em Combinatorial model categories have
presentations}, Adv. Math. {\bf 164}, no. 1 (2001), 177-201.

\bibitem[D4]{hend} D. Dugger, \emph{Spectral enrichments of model categories}, 
Homology Homotopy Appl. {\bf 8} (1), 2006.

\bibitem[DS1]{dereq} D. Dugger and B. Shipley, {\em $K$-theory and derived
equivalences},  Duke Math. J. {\bf 124} (2004), no. 3, 587--617.

\bibitem[DS2]{tpwe} D. Dugger and B. Shipley, {\em Topological
equivalences for differential graded algebras}, to appear in Advances
  in Math.

\bibitem[DS3]{stabmod} D. Dugger and B. Shipley,  {\em A curious
  example of two model categories and some associated differential
  graded algebras}, in preparation.

\bibitem[DK]{DK} W.G. Dwyer and D.M. Kan, {\em Function complexes in
homotopical algebra}, Topology {\bf 19} (1980), 427--440.

\bibitem[GS]{GS} J.~P.~C.~Greenlees and B.~Shipley, {\em An algebraic
model for rational torus-equivariant spectra}.  Preprint, 2004.
\verb!http://www.math.uic.edu/~bshipley/!

\bibitem[Hi]{H} P. Hirschhorn, {\em Model Categories and Their
Localizations}, Mathematical Surveys and Monographs, vol. 99,
Amer. Math. Soc., 2003.

\bibitem[Ho1]{hovey}
M.~Hovey, {\em Model categories}, Mathematical Surveys and Monographs,
{\bf 63}, American Mathematical Society,
Providence, RI, 1999, xii+209 pp.

\bibitem[Ho2]{stab} M.~Hovey, {\em Spectra and symmetric spectra in
general model categories\/}, J. Pure Appl. Algebra {\bf 165} (2001),
no. 1, 63--127.

\bibitem[HSS]{hss}
M.~Hovey, B.~Shipley, and J.~Smith, {\em Symmetric spectra},
J.\ Amer.\ Math.\ Soc.\ {\bf 13} (2000) 149--208.


\bibitem[ML]{MacL} S. Mac\,Lane, {\em Categories for the working
mathematician}, Graduate Texts in Math. {\bf 5}, Springer, New
York-Berlin, 1971.


\bibitem[SS1]{SS1}
S.~Schwede and B.~Shipley, {\em Algebras and modules
in monoidal model categories}, Proc.\ London Math.\ Soc.\ {\bf 80} (2000)
491-511.
                                                                               
\bibitem[SS2]{SS2}
S.~Schwede and B.~Shipley, {\em
Stable model categories are categories of modules},  Topology {\bf 42}
(2003), 103-153.

\bibitem[SS3]{SS3} S. Schwede and B. Shipley, {\em Equivalences of
monoidal model categories}, Algebr. Geom. Topol. {\bf 3} (2003),
287--334.


\bibitem[S]{S} B. Shipley, {\em $H\Z$-algebra spectra are
differential graded algebras}, to appear in Amer. Jour. Math.


\bibitem[T]{T} B. To\"en, {\em The homotopy theory of
dg-categories and derived Morita theory\/}, preprint, 2005.

\end{thebibliography}

\end{document}